\documentclass[smallextended]{svjour3}       % onecolumn (second format)
\smartqed  % flush right qed marks, e.g. at end of proof
\usepackage{graphicx}
\usepackage[T1]{fontenc}
\usepackage[utf8]{inputenc}
\usepackage{xcolor}
\usepackage{booktabs}
\usepackage{afterpage}
\usepackage{hyperref}
\usepackage{soul}

\usepackage{amsmath,amssymb}
\usepackage[linesnumbered,ruled,vlined]{algorithm2e}

\usepackage{newtxtext,newtxmath}

\usepackage{tikz}
\usepackage{pgfplots}
\usepgfplotslibrary{groupplots}
\newcommand{\dpair}[2]{\left\langle#1\right\rangle_{#2^{\ast}, #2}}
\newcommand{\iprod}[2]{\left(#1\right)_{#2}}
\newcommand{\norm}[2]{\left\lVert#1\right\rVert_{#2}}
\newcommand{\abs}[1]{\left\lvert#1\right\rvert}
\newcommand{\ud}{\mathrm{d}}
\newcommand{\eye}{\mathrm{I}}

\DeclareMathOperator{\kernel}{ker}
\DeclareMathOperator{\range}{ran}

\newtheorem{assumption}{Assumption}

\SetKwFor{Loop}{Loop}{}{EndLoop}
\SetKw{Break}{break}

\journalname{Math. Prog., Ser. A}

\smartqed

\begin{document}

\title{A Sequential Homotopy Method for Mathematical Programming Problems%
\thanks{We gratefully acknowledge support by the German Federal
Ministry of Education and Research under grants MOPhaPro (05M16VHA) and MOReNet
(05M18VHA). HGB also gratefully acknowledges support by the German
Research Foundation within the priority program DFG-SPP 1962. We thank Paul
Manns for his helpful comments on an earlier version of the manuscript.}
}


\author{Andreas Potschka \and Hans Georg Bock}


\institute{A. Potschka, H.G. Bock \at
              Im Neuenheimer Feld 205, 69120 Heidelberg, Germany\\
              Tel.: +49-6221-5414637\\
              Fax: +49-6221-5414629\\
              \email{\{potschka|bock\}@iwr.uni-heidelberg.de}           %  \\
}

\date{Received: date / Accepted: date}

\maketitle

\begin{abstract}
  We propose a sequential homotopy method for the solution of mathematical
  programming problems formulated in abstract Hilbert spaces under the Guignard
  constraint qualification. The method is equivalent to performing projected
  backward Euler timestepping on a projected gradient/antigradient flow of the
  augmented Lagrangian. The projected backward Euler equations can be
  interpreted as the necessary optimality conditions of a primal-dual proximal
  regularization of the original problem. The regularized problems are always
  feasible, satisfy a strong constraint qualification guaranteeing uniqueness of
  Lagrange multipliers, yield unique primal solutions provided that the
  stepsize is sufficiently small, and can be solved by a continuation in the
  stepsize. We show that equilibria of the projected gradient/antigradient flow
  and critical points of the optimization problem are identical, provide
  sufficient conditions for the existence of global flow solutions, and show
  that critical points with emanating descent curves cannot be asymptotically
  stable equilibria of the projected gradient/antigradient flow, practically
  eradicating convergence to saddle points and maxima. The sequential homotopy
  method can be used to globalize any locally convergent optimization method
  that can be used in a homotopy framework. We demonstrate its efficiency for a
  class of highly nonlinear and badly conditioned control constrained elliptic
  optimal control problems with a semismooth Newton approach for the regularized
  subproblems.

\end{abstract}

\section{Introduction}
\label{sec:intro}

Let $X$ and $Y$ be real Hilbert spaces and $C \subseteq X$ a nonempty closed
convex set. Let the nonlinear objective function $\phi: X \to \mathbb{R}$ and
the nonlinear constraint function $c: X \to Y$ be twice continuously
Fr\'echet differentiable. We consider the mathematical programming problem
\begin{equation}
  \min \phi(x) \quad \text{over } x \in C \quad \text{subject to } c(x) = 0.
  \label{eqn:mp}
\end{equation}
This formulation is equivalent to a more prevalent formulation that allows $c(x)
\in C_c$ for some nonempty closed convex set $C_c$ (by the use of slack
variables $s \in Y$ via $c(x) - s = 0$ and $(x, s) \in C \times C_c$). Further
restrictions on the overall setting are stated in section~\ref{sec:assumptions}
after we settle the notation in section~\ref{sec:notation}.

This setting naturally comprises finite dimensional problems (also known as
Nonlinear Programming Problems, NLPs) of the form
\begin{equation*}
  \min_{x \in \mathbb{R}^n} \phi(x) \quad \text{subject to } x^\mathrm{l} \le
  x \le x^\mathrm{u} \text{ and } c(x) = 0
\end{equation*}
with $X = \mathbb{R}^n$, $Y = \mathbb{R}^m$, and $C = \{ x \in \mathbb{R}^n
\mid x^\mathrm{l} \le x \le x^\mathrm{u} \}$, where some components of
$x^\mathrm{u}$ and $x^\mathrm{l}$ may take on values of $\pm \infty$.

Another popular example is partial differential equation (PDE) constrained
optimization, where $X = U \times Q$ is a product of the state and control
space, $C$ encodes pointwise constraints on the controls, and $c(x) = c((u, q))
= 0$ is the PDE constraint, where we often assume that the state $u \in U$ is
locally uniquely determined by the control $q \in Q$ as an implicit function
$u(q)$ via $c((u(q), q)) = 0$.

\subsection{Structure of the article}

We give a concise overview of the results of this article in
section~\ref{sec:overview}.
We outline our contributions and connections to existing methods
in section~\ref{sec:contributions}. In the remainder of section~\ref{sec:intro},
we settle our notation, state the general assumptions, and provide the
statements of important classical results. We give a short proof of the
necessary optimality conditions we use and discuss two central constraint
qualifications in section~\ref{sec:noc}. Our main results on projected
gradient/antigradient flows for~\eqref{eqn:mp} follow in
section~\ref{sec:gradflow}. The application of a projected backward Euler method
on the projected gradient/antigradient flow results in a sequential homotopy
method, which we describe in section~\ref{sec:seqhom}. We present numerical
results for a local semismooth Newton method globalized by the sequential
homotopy approach for a class of highly nonlinear and badly conditioned elliptic
PDE-constrained optimal control problems with control constraints in
section~\ref{sec:PDE_opt}.

\subsection{Overview: A novel solution approach based on a sequence of
homotopies}
\label{sec:overview}

We propose the following general solution approach in this paper: We construct
and analyze existence and uniqueness of a primal-dual projected
gradient/antigradient flow for an augmented Lagrangian.
The equilibria of the flow are critical points of~\eqref{eqn:mp}
and vice versa.  Under reasonable assumptions, we prove that critical points
that are not local minima cannot be asymptotically stable. Small perturbations
will make the flow escape these unwanted critical points. We then apply a
projected version of backward Euler timestepping. We provide an interpretation
of the backward Euler equations as the optimality conditions of a primal-dual
proximally regularized counterpart to~\eqref{eqn:mp}, which satisfies a strong
constraint qualification, even though~\eqref{eqn:mp} might only satisfy the
Guignard constraint qualification~\cite{Guignard1969}, the weakest of all
constraint qualifications.  This gives rise to a sequential homotopy method, in
which a sequence of proximally regularized subproblems needs to be solved by
(possibly inexact) fast numerical methods that are only required to converge
locally.

%%% Arxiv version (including supplementary material):
%Before we describe this apporach in full detail and rigor, we would like to
%sketch without proofs its salient features with an illustrative example in
%finite dimensions without inequalities. We assume the reader is familiar with
%standard NLP theory.
%
%\input{pendulum}

%% Journal version (excluding supplementary material)
We invite the reader to read the supplementary material, in which we
sketch without proofs the salient features of our approach with an illustrative
example in finite dimensions without inequalities.

\subsection{Related work and contributions}
\label{sec:contributions}

We advance and bridge several fields of optimization with this paper.

The field of globalized Newton methods based on differential equation methods
applied to the Newton flow started in the early 1950s with
Davidenko~\cite{Davidenko1953} and continues to raise scientific interest over
the
decades~\cite{Deuflhard1974,Ascher1987,Hohmann1994,Deuflhard1998,Bock2000,Potschka2016,Lubkoll2017,Potschka2018,Deuflhard2018},
predominantly due to the affine invariance properties of the Newton
flow~\cite{Deuflhard2004}. 
By trading the affine invariance of the Newton flow for the stability properties
of the gradient flow, we obtain from a dynamical systems point of view the
advantage of being repelled from maxima or saddle points when solving nonlinear
optimization problems.

The Newton method, which is
equivalent to forward Euler timestepping on the Newton flow with stepsize
$\Delta t = 1$, has the prominent property of quadratic local convergence.
Backward Euler timestepping on the gradient/antigradient flow can attain
superlinear local convergence if the solution is sufficiently regular so that we
can take the stepsize $\Delta t$ to infinity or, equivalently, drive the
proximal coefficient $\lambda$ to zero, provided that we use a local solver in
the numerical homotopy method with at least superlinear local convergence.
Driving $\lambda$ to zero is usually possible if the solution satisfies certain
second order sufficient optimality conditions.

Three methods in the field of convex optimization are closely related to our
approach.  The first method is the proximal point algorithm for closed proper
convex functions, which can be interpreted as a backward Euler timestepping on
the gradient flow of the objective function, while the gradient descent method
amounts to forward Euler timestepping on the gradient flow (see, e.g.,
\cite[sec.~4.1]{Parikh2014} and references therein). We extend this approach to
nonconvex optimization problems with explicit handling of nonlinear equality
constraints, as they appear for instance in optimal control. To this end, we extend a
second method, the primal-dual projected gradient/antigradient flow of
\cite[chap.~6, 7]{Arrow1958}, from the finite-dimensional convex to the
infinite-dimensional nonconvex setting with the help of an augmented Lagrangian
technique in the framework of projected differential equations in Hilbert
space~\cite{Cojocaru2004}.
The third method we extend is the closely related Arrow--Hurwicz gradient method
\cite[chap.~10]{Arrow1958}, which amounts to projected forward Euler
timestepping on the projected gradient/antigradient flow of the Lagrangian
without augmentation ($\rho = 0$). Our sequential homotopy method is equivalent
to projected backward Euler timestepping. Hence, it bears the same connection
with the Arrow--Hurwicz gradient method as the proximal point algorithm with
gradient descent.

From a Sequential Quadratic Programming (SQP) perspective (see,
e.g.,~\cite{Nocedal2006}), our approach resolves all the numerical difficulties
on the nonlinear level such as subproblem infeasibility, degeneracy, and
nonconvexity due to indefinite subproblem Hessians.
Existing approaches often pass these difficulties on to
the level of the quadratic subproblem solvers, which may fail to resolve these
issues in a way that guarantees convergence of the overall nonlinear iteration.
Our method can thus be used as a black-box globalization framework for any
locally convergent optimization method that can be used within a continuation
framework, e.g., methods of structure-exploiting inexact Sequential Quadratic
Programming
(SQP)~\cite{Hintermueller2004,Potschka2012,Potschka2013,Potschka2015,Hante2015}
or semismooth Newton
methods~\cite{Mifflin1977,Qi1993,Ulbrich2002,Hintermueller2002,Ito2004,Hintermueller2004,Ulbrich2011,Hintermueller2010}.
The local methods are even allowed to converge to maxima and saddle points.
These issues are taken care of by our sequential homotopy method. For the
application of local SQP methods, we can guarantee that the quadratic
subproblems are always feasible and that they satisfy a strong constraint
qualification that implies unique subproblem Lagrange multipliers. In addition,
they are convex if the augmentation parameter $\rho$ is sufficiently large and
the stepsize $\Delta t$ is sufficiently small when we are still far away from a
solution.

Our approach uses the theory of projected differential equations due to Cojocaru
and Jonker~\cite{Cojocaru2004}, which have a tight connection to differential
inclusions~\cite{Aubin1984} and evolutionary/differential variational
inequalities~\cite{Cojocaru2005,Pang2008}. We are mainly interested in their
equilibrium points, which satisfy a variational inequality (VI). Other methods
to compute solutions to VIs have been described in the literature (see,
e.g.,~\cite{Munson2001,Benson2006}), which are based on semismooth
iterations on reformulations using special Nonlinear Complementarity Problem
(NCP) functions.

Projected gradient flows for constrained optimization problems in finite
dimensions have also been considered with techniques from Riemannian geometry
(see, e.g.,~\cite{Jongen2001,Shikhman2009,Hauswirth2016,Hauswirth2018} and
references therein), but the
resulting methods produce only feasible iterates. It is often computationally
wasteful to satisfy all constraints for iterates far away from an optimum and to
force the iterates to follow a feasible manifold with possibly high curvature.

For an introduction to augmented Lagrangian approaches in Hilbert spaces we
refer to~\cite{Ito2008} and references therein. We point out that our approach
relies on the augmented Lagrangian mainly to remove negative curvature of the
Lagrangian in the kernel of the constraints. In contrast to classical augmented
Lagrangian methods, we do not alternate between updates of the primal and dual
variables but rather update primal and dual variables simultaneously as in
augmented Lagrangian-SQP methods~\cite[chap.~6]{Ito2008}.

\subsection{Notation}
\label{sec:notation}
We abbreviate the nonnegative real numbers with $\mathbb{R}_{\ge 0}$.
By $(x_k) \subset X$ we denote a sequence $x_0, x_1, \dotsc$ of elements in $X$.
By $X^{\ast}$ we denote the topological dual of $X$,
by $\iprod{.,.}{X}: X \times X \to \mathbb{R}$ the inner product,
by $\norm{.}{X}: X \to \mathbb{R}_{\ge 0}$ the norm, and
by $\dpair{.,.}{X}: X^{\ast} \times X \to \mathbb{R}$ the duality pairing.
By $R_X: X^{\ast} \to X$ we denote the Riesz isomorphism (see, e.g.,
\cite[sec.~III.6]{Yosida1995}), which satisfies the identity
\begin{align*}
  \iprod{R_X x^\ast, x}{X} &= \dpair{x^{\ast}, x}{X} \quad
  \text{for all } x^{\ast} \in X^{\ast}, x \in X
\end{align*}
and likewise for $Y$.
As usual, $\mathcal{L}(X, Y)$ denotes the Banach-space of all continuous
linear operators from $X$ to $Y$.
For $A \in \mathcal{L}(X, Y)$, the (Banach space) dual operator $A^{\ast} \in
\mathcal{L}(Y^{\ast}, X^{\ast})$ and the (Hilbert space) adjoint operator
$A^{\star} \in \mathcal{L}(Y, X)$ are defined by
\begin{align*}
  \dpair{A^{\ast} y^{\ast}, x}{X} &= \dpair{y^{\ast}, A x}{Y}
  &&\text{for all } x \in X, y^{\ast} \in Y^{\ast},\\
  \iprod{A^{\star} y, x}{X} &= \iprod{y, A x}{Y}
  &&\text{for all } x \in X, y \in Y,
\end{align*}
which implies $A^{\star} R_{Y} = R_{X} A^{\ast}$.
We denote the Fr\'echet-derivative of $c(x)$ with $c'(x) \in \mathcal{L}(X,
Y)$.  We denote the objective gradient by $\nabla \phi(x) = R_{X} \phi'(x) \in
X$ and the adjoint of the constraint derivative by $\nabla c(x) = \left( c'(x)
\right)^{\star} \in \mathcal{L}(Y, X)$.
For a linear operator $A \in \mathcal{L}(X, Y)$, we denote its kernel by
$\kernel(A) = \{x \in X \mid A x = 0 \}$ and its range by $\range(A) = \{ y \in
Y \mid \exists x \in X: y = A x \}$.
%The abbreviation $\cvx(M)$ denotes the the convex hull of a set $M$ and
%$\ccvx(M)$ the closure of $\cvx(M)$ within its corresponding space, which will
%be clear from the context.
For an open set $\Omega \in \mathbb{R}^n$, we denote with $L^2(\Omega)$ the
standard Hilbert space of square Lebesgue-integrable functions on $\Omega$, with
$H^1_0(\Omega)$ the Sobolev-space of functions with square Lebesgue-integrable
derivatives and zero trace at the boundary, and with $H^{-1}(\Omega)$ its dual
space.
We denote the feasible set of~\eqref{eqn:mp} with $\mathcal{F} = \{ x \in C \mid c(x) = 0 \}$.

\subsection{General assumptions}
\label{sec:assumptions}
A central role in this article is played by the augmented objective and
augmented Lagrangian
\begin{align}
  \label{eqn:aug_Lagrangian}
  \phi^{\rho}(x) &= \phi(x) + \frac{\rho}{2} \norm{c(x)}{Y}^2, &
  L^{\rho}(x,y) &= \phi^{\rho}(x) + \iprod{y, c(x)}{Y},
\end{align}
defined for some fixed $\rho \in \mathbb{R}_{\ge 0}$ and arbitrary $x \in C$ and
$y \in Y$.
Throughout this article, we make the following assumptions: 
\begin{assumption}
  \label{ass:closed_range}
  For all $x \in \mathcal{F}$, $\range(c'(x))$ is closed in $Y$.
\end{assumption}
\begin{assumption}
  \label{ass:phi_rho_bounded}
  For some fixed $\rho \in \mathbb{R}_{\ge 0}$ we have the coercivity condition
  \[
    \phi^{\rho}_{\mathrm{low}} = \inf_{x \in C} \phi^{\rho}(x) > -\infty
    \quad \text{and} \quad
    \lim_{\norm{x}{X} \to \infty} \phi^{\rho}(x) = \infty.
  \]
\end{assumption}
\begin{assumption}
  \label{ass:lipschitz}
  The functions $c(x)$, $L^{\rho}(x,y)$ and the gradient $\nabla L^{\rho}(x,y)$
  are locally Lipschitz continuous.
\end{assumption}

\subsection{Well-known results}
Let us recall the following well-known definitions.
\begin{definition}[Tangent cone]
  For $\bar{x} \in X$ and a nonempty set $M \subseteq X$, we call
  \begin{align*}
    T(M, \bar{x}) = \{ d \in X \mid~&\text{there exist sequences } (x_k) \subset
    M, (\lambda_k) \subset \mathbb{R}_{\ge 0}\\
    &\text{with } x_k \to \bar{x} \text{ and }
    \lambda_k (x_k - \bar{x}) \to d \text{ as } k \to \infty \}
  \end{align*}
  the \emph{tangent cone} to $M$ at $\bar{x}$.
\end{definition}
\begin{definition}[Projection]
  For a nonempty closed convex set $K \subseteq X$,
  we denote by $P_{K}: X \to K$ the \emph{projection operator} of $X$ onto $K$,
  which is uniquely defined by
  \[
    \norm{P_{K}(x) - x}{X} = \inf_{\tilde{x} \in K} \norm{\tilde{x} - x}{X}
    \quad \text{for all } x \in X.
  \]
\end{definition}
For properties of projection operators, we refer the reader to
\cite{Zarantonello1971}.
\begin{definition}[Polar cone]
  For a cone $K \subseteq X$, we call
  \[
    K^{-} = \left\{ d \in X \mid \iprod{d, x}{X} \le 0
    \text{ for all } x \in K \right\}
  \]
  the \emph{polar cone} of $K$.
\end{definition}

\begin{remark}
  If $K \subseteq X$ is a linear subspace, then $x \in K$ implies $-x \in K$ and
  thus equality holds in the definition of $K^{-} = \left\{ d \in
  X \mid \iprod{d, x}{X} = 0 \text{ for all } x \in K \right\} = K^{\perp}.$
\end{remark}

We shall make use of the following classical results from convex analysis.
\begin{lemma}[Moreau decomposition]
  \label{lem:moreau_decomp}
  If $K \subseteq X$ is a nonempty closed convex cone, then every $x \in X$ has
  a unique decomposition $x = P_{K}(x) + P_{K^{-}}(x) =: x^{+} +
  x^{-}$, where $\iprod{x^{-}, x^{+}}{X} = 0$. A simple consequence is the
  identity 
  \[
    \iprod{x, P_{K}(x)}{X} = \iprod{x^{+} + x^{-}, x^{+}}{X}
    = \norm{P_{K}(x)}{X}^2.
  \]
\end{lemma}
\begin{proof}
  See \cite{Moreau1962} according to \cite[Lemma~2.2 and
  Corollary~2]{Zarantonello1971}.
  \qed
\end{proof}

\begin{lemma}
  \label{lem:polar_cone_and_projection}
  Let $K \subseteq X$ be a nonempty closed convex set and let $\bar{x} \in
  K$. If $x \in T^{-}(K, \bar{x}) + \bar{x}$, then $P_{K}(x) = \bar{x}$.
\end{lemma}
\begin{proof}
  Choose any $y \in K$. Then, $y - \bar{x} \in T(K, \bar{x})$, e.g., with
  $\lambda_{k} = k + 1$ and $x_{k} = (1-\lambda_{k}^{-1}) \bar{x} +
  \lambda_{k}^{-1} y \in K$. Because $x - \bar{x} \in T^{-}(K, \bar{x})$, we
  obtain $\iprod{x - \bar{x}, y - \bar{x}}{X} \le 0.$ The result follows
  from~\cite[Lemma~1.1]{Zarantonello1971}, because $y \in K$ was chosen
  arbitrarily.
  \qed
\end{proof}

\section{Necessary optimality conditions}
\label{sec:noc}

The basis for the sequential homotopy method we propose in Sec.~\ref{sec:seqhom}
is a necessary optimality condition due to Guignard \cite{Guignard1969}. Because
the separation of nonlinearities $c(x) = 0$ and inequalities $x \in C$
in~\eqref{eqn:mp} allow for a much shorter proof, we state it here for the sake
of convenience.

%\begin{definition}[Pseudotangent cone]
%  For $\bar{x} \in X$ and a nonempty set $M \subseteq X$, we call $P(M,\bar{x})
%  = \ccvx(T(M, \bar{x}))$ the \emph{pseudotangent cone} to $M$ at $\bar{x}$.
%\end{definition}

\begin{lemma}
  \label{lem:gen_noc}
  If $\bar{x} \in \mathcal{F}$ is a local optimum of \eqref{eqn:mp}, then
  $-\nabla \phi(\bar{x}) \in T^{-}(\mathcal{F}, \bar{x})$.
\end{lemma}
\begin{proof}
  Let $d \in T(\mathcal{F}, \bar{x})$ with corresponding sequences $(x_k)
  \subset X$ and $(\lambda_k) \subset \mathbb{R}_{\ge 0}$. Using the shorthand
  $d_k = \lambda_k (x_k - \bar{x})$, we obtain the assertion from letting $k \to
  \infty$ in
  \[
    0 \le \lambda_k \left[ \phi(x_k) - \phi(\bar{x}) \right]
    = \dpair{\phi'(\bar{x}), d_k}{X} + \norm{d_k}{X} \frac{o\left(\norm{x_k -
    \bar{x}}{X}\right)}{\norm{x_k - \bar{x}}{X}}
    \to \iprod{\nabla \phi(\bar{x}), d}{X}.~
    \qed
  \]
\end{proof}

\begin{definition}[GCQ]
  \label{def:GCQ}
  We say that the \emph{Guignard Constraint Qualification (GCQ)} holds at
  $\bar{x} \in \mathcal{F}$ if 
  \[
    \kernel^{\perp} (c'(\bar{x})) + T^{-} (C, \bar{x}) = T^{-}(\mathcal{F},
    \bar{x}).
  \]
\end{definition}

\begin{theorem}[Necessary optimality conditions]
  \label{thm:noc}
  If $\bar{x} \in \mathcal{F}$ is a local optimum of \eqref{eqn:mp} that
  satisfies GCQ, then there exists a multiplier $\bar{y} \in Y$ such that
  \begin{equation}
    -\nabla \phi(\bar{x}) - \nabla c(\bar{x}) \bar{y} \in T^{-}(C, \bar{x}).
    \label{eqn:critical}
  \end{equation}
\end{theorem}
\begin{proof}
  The proof is based on the Closed Range Theorem
  (see, e.g., \cite[sec.~VII.5]{Yosida1995} with premultiplication by the Riesz
  isomorphism $R_X$ to obtain the Hilbert space version), which states that
  Assumption~\ref{ass:closed_range} is equivalent to
  \[
    \kernel^{\perp}(c'(\bar{x})) = \range(\nabla c(\bar{x})).
  \]
  Together with Lemma~\ref{lem:gen_noc} and GCQ we obtain
  \[
    -\nabla \phi(\bar{x}) \in T^{-}(\mathcal{F}, \bar{x})
    = \kernel^{\perp}(c'(\bar{x})) + T^{-}(C, \bar{x})
    = \range(\nabla c(\bar{x})) + T^{-}(C, \bar{x}).
  \]
  Thus, there exists a $\bar{y} \in Y$ such that
  $-\nabla \phi(\bar{x}) - \nabla c(\bar{x}) \bar{y} \in T^{-}(C, \bar{x}).$
  \qed
\end{proof}

\begin{definition}[Critical point]
  We call $(\bar{x}, \bar{y}) \in \mathcal{F} \times Y$ a \emph{critical point}
  if \eqref{eqn:critical} holds.
\end{definition}

The method we propose below enjoys the benefit that its subproblems lift the
original problem into a larger space with additional structural properties in
$X$, $C$, and $c$, which result in satisfaction of a constraint qualification
that is much stronger than GCQ, even though problem~\eqref{eqn:mp} only
satisfies GCQ.
\begin{lemma}%{Implicit Function Constraint Qualification}
  \label{lem:IFCQ}
  Let $X = U \times Q$, equipped with the canonical inner product derived from
  the Hilbert spaces $U$ and $Q$, and let $C = U \times C_{Q}$ for some nonempty
  closed convex set $C_{Q} \subseteq Q$. Furthermore, assume there exists a
  continuously Fr\'echet-differentiable mapping $S: C_{Q} \to U$ such that
  for all $x = (u, q) \in C$
  \begin{align*}
    \text{\rm (a)} ~&~ c( (u, q)) = 0 ~ \text{ iff } ~ u = S(q), &
    \text{\rm (b)} ~&~ \range c'_{u}(x) = Y, &
    \text{\rm (c)} ~&~ \range \nabla_{u} c(x) = U.
  \end{align*}
  Then, $\mathcal{F}$ is nonempty, every $\bar{x} \in \mathcal{F}$ satisfies
  GCQ, and the Lagrange multiplier $\bar{y}$ in~\eqref{eqn:critical} is uniquely
  determined.
\end{lemma}
\begin{proof}
  The feasible set $\mathcal{F} = \{ (S(q), q) \mid q \in C_Q \}$ is nonempty
  because $C_Q$ is nonempty. Let $\bar{x} = (S(\bar{q}), \bar{q}) \in
  \mathcal{F}$ and choose some $d \in T(C_{Q}, \bar{q})$. By definition, there
  exist sequences $(q_k) \subset C_{Q}$ and $(\lambda_k) \subset \mathbb{R}_{\ge
  0}$ such that $\lambda_k (q_k - \bar{q}) \to d$. Using (a), we choose a
  sequence $(x_k) \subset \mathcal{F}$ according to $x_k = (S(q_k), q_k)$ to
  guarantee $x_k \to \bar{x}$ and
  \begin{equation}
    \label{eqn:tangent_sequence}
    \begin{aligned}
      \lambda_k (x_k - \bar{x})
      &= \lambda_k (S(q_k) - S(\bar{q}), q_k - \bar{q})\\
      &= \lambda_k (S'(\bar{q}) (q_k - \bar{q}) + o(\norm{q_k - \bar{q}}{Q}),
      q_k - \bar{q}) \to (S'(\bar{q}) d, d),
    \end{aligned}
  \end{equation}
  which shows that $T(\mathcal{F}, \bar{x}) \supseteq \{ (S'(\bar{q}) d, d) \mid
  d \in T(C_{Q}, \bar{q}) \}.$
  In order to show that equality holds between the two sets, we notice that if
  $(e, d) \in T(\mathcal{F}, \bar{x})$ then $d \in T(C_{Q}, \bar{q})$
  and~\eqref{eqn:tangent_sequence} implies $e = S'(\bar{q}) d$. Hence, we obtain
  \begin{equation*}
    T(\mathcal{F}, \bar{x}) = \{ (S'(\bar{q}) d, d) \mid d \in
    T(C_{Q}, \bar{q}) \}.
  \end{equation*}
  In order to compute its polar cone, let $x = (u, \tilde{q}) \in X$ such that
  \[
    0 \ge \iprod{u, S'(\bar{q}) d}{U} + \iprod{\tilde{q}, d}{Q}
    = \iprod{\nabla S(\bar{q}) u + \tilde{q}, d}{Q} \quad
    \text{for all } d \in T(C_{Q}, \bar{q}).
  \]
  We choose $q = \nabla{S}(\bar{q}) u + \tilde{q}$ in order to obtain
  \begin{equation}
    \label{eqn:IFCQ_polar_tangent_cone}
    \begin{aligned}
      T^{-}(\mathcal{F}, \bar{x})
      &= \left\{ (u, \tilde{q}) \in X \mid \iprod{u, e}{U} 
      + \iprod{\tilde{q},d}{Q} \le 0
      \text{ for all } (e,d) \in T(\mathcal{F}, \bar{x}) \right\}\\
      &= \left\{ (u, \tilde{q}) \in X \mid
      \iprod{\nabla S(\bar{q}) u + \tilde{q}, d}{Q} \le 0
      \text{ for all } d \in T(C_{Q}, \bar{q}) \right\}\\
      &= \left\{ (u, q - \nabla S(\bar{q}) u) \mid u \in U, q \in
      T^{-}(C_{Q}, \bar{q}) \right\}.
    \end{aligned}
  \end{equation}
  For the other polar cone in the definition of GCQ, we get
  \begin{equation}
    T^{-}(C, \bar{x}) = T^{-}(U \times C_{Q}, (\bar{u}, \bar{q}))
    = \left(U \times T(C_{Q}, \bar{q}) \right)^{-}
    = \{ 0 \} \times T^{-}(C_{Q}, \bar{q}).
    \label{eqn:IFCQ_polar_cone}
  \end{equation}
  Taking the derivative of $c(S(q), q) = 0$ with respect to $q$ in
  direction $d \in Q$ yields
  \[
    c'_{u}(\bar{x}) S'(\bar{q}) d + c'_{q}(\bar{x}) d = 0.
  \]
  As a consequence of the Closed Range Theorem~\cite[sec.~VII.5,
  Corollary~1]{Yosida1995}, (c) is equivalent to the existence of a continuous
  inverse of $c'_{u}(\bar{x})$, from which we see that
  \[
    \kernel c'(\bar{x}) = 
    \left\{ (e,d) \in X \mid c'_{u}(\bar{x}) e + c'_{q}(\bar{x}) d = 0 \right\}
    = \left\{ (S'(\bar{q}) d, d) \mid d \in Q \right\}.
  \]
  Thus, its orthogonal complement amounts to
  \begin{equation}
    \begin{aligned}
      \kernel^{\perp} c'(\bar{x})
      &= \left\{ (u,q) \in X \mid \iprod{u, S'(\bar{q}) d}{U} + \iprod{q,d}{Q} =
      0 \text{ for all } d \in Q \right\}\\
      &= \left\{ (u,q) \in X \mid \iprod{\nabla S(\bar{q}) u + q, d}{Q} = 0
      \text{ for all } d \in Q \right\}\\
      &= \left\{ (u, -\nabla S(\bar{q}) u) \mid u \in U \right\}.
    \end{aligned}
    \label{eqn:IFCQ_kernel_complement}
  \end{equation}
  Hence, it follows from \eqref{eqn:IFCQ_polar_cone},
  \eqref{eqn:IFCQ_kernel_complement}, and \eqref{eqn:IFCQ_polar_tangent_cone}
  that
  \[
    T^{-}(C, \bar{x}) + \kernel^{\perp} c'(\bar{x}) = 
    \left\{ (u, q - \nabla S(\bar{q}) u \mid u \in U, q \in T^{-}(C_{Q},
    \bar{q}) \right\}
    = T^{-}(\mathcal{F}, \bar{x}),
  \]
  which shows that GCQ holds at $\bar{x}$. Regarding multiplier uniqueness, we
  take the $U$-components of~\eqref{eqn:critical}
  and~\eqref{eqn:IFCQ_polar_cone} to deduce
  \[
    \nabla_{u} \phi(\bar{x}) + \nabla_{u} c(\bar{x}) \bar{y} = 0,
  \]
  from which the uniqueness of $\bar{y}$ follows from the the existence of a
  continuous inverse of $\nabla_{u} c(\bar{x})$ by virtue of (b)
  and~\cite[sec.~VII.5, Corollary~1]{Yosida1995}.
  \qed
\end{proof}

\section{Projected gradient/antigradient flow}
\label{sec:gradflow}

We study a primal-dual gradient/anti-gradient flow (from now on simply called
\emph{gradient flow}) of the augmented Lagrangian $L^{\rho}$, defined in
\eqref{eqn:aug_Lagrangian}, projected on the closed convex set $C$ in the
framework of projected differential equations in Hilbert space
\cite{Cojocaru2004} according to
\begin{align}
  \label{eqn:gradflow}
  \dot{x}(t) &= P_{T(C,x(t))} \left(-\nabla_x L^{\rho}(x(t), y(t))
  \right), &
  \dot{y}(t) &= \nabla_y L^{\rho}(x(t), y(t)),
\end{align}
where the gradients with respect to $x$ and $y$ evaluate to
\begin{align*}
  \nabla_{x} L^{\rho}(x, y) &= \nabla \phi(x) + \nabla c(x) \left[ y + \rho c(x)
  \right], &
  \nabla_{y} L^{\rho}(x, y) &= c(x).
\end{align*}
The following existence theorem uses $L^{\rho}$ and $\frac{1}{2}
\norm{c(.)}{Y}^2$ as Lyapunov-type functions.
Due to Lemma~\ref{lem:moreau_decomp}, the $t$-derivative of $L^{\rho}$
along the flow is given by
\begin{equation}
  \label{eqn:dLdt}
  \begin{aligned}
    \frac{\ud}{\ud t} L^{\rho}(x(t), y(t))
    &= \iprod{\nabla_x L^{\rho}(x(t), y(t)), \dot{x}(t)}{X}
    + \iprod{\nabla_y L^{\rho}(x(t), y(t)), \dot{y}(t)}{Y}\\
    &= -\norm{P_{T(C, x(t))}\left(-\nabla_x L^{\rho}(x(t), y(t))\right)}{X}^2
    + \norm{c(x(t))}{Y}^2.
  \end{aligned}
\end{equation}
The positive sign in front of the last term in~\eqref{eqn:dLdt} reflects
the saddle point nature of the Lagrangian approach and complicates the use of
Lyapunov arguments in comparison to the unconstrained case. We pursue the basic
idea that by increasing $\rho$, we can make the negative term overpower the
$\rho$-independent positive term. That this is not always possible will be
discussed after the following theorem. 

\begin{theorem}[Unique existence of solutions]
  \label{thm:existence_of_trajectories}
  Let Assumptions~\ref{ass:phi_rho_bounded}
  and~\ref{ass:lipschitz} be satisfied.
  Then, there exists an interval $[0, t_{\mathrm{final}}]$ and a uniquely
  determined pair of absolutely continuous functions $(x, y): [0,
  t_{\mathrm{final}}] \to C \times Y$
  that satisfy the projected gradient flow equation \eqref{eqn:gradflow} and
  $(x(0), y(0)) = (x_0, y_0)$. The final time $t_{\mathrm{final}}$ can be
  extended as long as the condition
  \begin{align}
    \label{eqn:L_descent}
    \frac{\ud}{\ud t} L^{\rho}(x(t), y(t))
    &\le 0
  \end{align}
  holds almost everywhere on $[0, t_{\mathrm{final}}]$. In addition,
  if for some $\gamma_1, \gamma_2 \in (0,1)$ the
  conditions~\eqref{eqn:L_descent} and
  \begin{align}
    \label{eqn:c_descent}
    \gamma_1 \frac{\ud}{\ud t} \left( \frac{1}{2} \norm{c(x(t))}{Y}^2 \right)
    &\le -\frac{\ud}{\ud t} L^{\rho}(x(t), y(t))
    -\gamma_2 \norm{c(x(t))}{Y}^2,
  \end{align}
  hold almost everywhere in $\mathbb{R}_{\ge 0}$, we have
  \begin{equation}
    \label{eqn:bounded_in_L2}
    \int_{0}^{\infty} \norm{P_{T(C, x(t))}\left( -\nabla_x L^{\rho}(x(t), y(t))
    \right)}{X}^2 \ud t < \infty
    \quad \text{and} \quad
    \int_{0}^{\infty} \norm{c(x(t))}{Y}^2 \ud t < \infty.
  \end{equation}
  Furthermore, if there is a set $M \subseteq X \times Y$ such that $\nabla
  L^{\rho}$ is (globally) Lipschitz continuous on $M$ and $(x(t), y(t)) \in M$
  for all $t \in [0, \infty)$, we obtain
  \begin{equation*}
    P_{T(C, x(t))}\left( -\nabla_x L^{\rho}(x(t), y(t)) \right) \to 0
    \quad \text{and} \quad
    c(x(t)) \to 0
    \quad \text{for } t \to \infty.
  \end{equation*}
\end{theorem}
\begin{proof}
  By Assumption~\ref{ass:lipschitz}, $\nabla L^{\rho}(x, y)$ is Lipschitz
  continuous in a neighborhood of $(x_0, y_0)$ with some Lipschitz constant $b <
  \infty$.
  By virtue of \cite[Theorem~3.1]{Cojocaru2004}, there exists an $l > 0$ and a
  uniquely determined pair of absolutely continuous functions $(x,y): [0, l] \to
  C \times Y$ that satisfy \eqref{eqn:gradflow} for almost all $t \in [0, l]$
  and $x(0) = x_0$, $y(0) = y_0$. 
  Without loss of generality, \eqref{eqn:L_descent} is satisfied on $[0, l]$ and
  we can repeatedly extend the local solution by the above arguments until
  \eqref{eqn:L_descent} or \eqref{eqn:c_descent} is violated for
  some $t_{\mathrm{final}} > 0$. As long as \eqref{eqn:L_descent} is satisfied,
  no blowup is possible in finite time. To see this, we first observe that
  \begin{align}
    \label{eqn:norm_y_derivative}
    \frac{\ud}{\ud t} \left( \frac{1}{2} \norm{y(t)}{Y}^2 \right)
    &= \iprod{y(t), c(x(t))}{Y},
  \end{align}
  which implies in combination with \eqref{eqn:L_descent} and
  Assumption~\ref{ass:phi_rho_bounded} that
  \begin{align*}
    \frac{1}{2} \norm{y(t)}{Y}^2 &= \frac{1}{2} \norm{y_0}{Y}^2
    + \int_{0}^{t} \iprod{y(\tau), c(x(\tau))}{Y} \ud \tau\\
    &= \frac{1}{2} \norm{y_0}{Y}^2 + \int_{0}^{t} \left[ L^{\rho}(x(\tau),
    y(\tau)) - \phi^{\rho}(x(\tau)) \right] \ud \tau\\
    &\le \frac{1}{2} \norm{y_0}{Y}^2 + t \left[ L^{\rho}(x_0, y_0) -
    \phi^{\rho}_{\mathrm{low}} \right].
  \end{align*}
  This establishes that there can be no blowup of $y$ in finite time. In
  addition, $x$ cannot blow up in finite time because then $L^{\rho}(x(t),
  y(t))$ would tend to infinity by virtue of
  Assumption~\ref{ass:phi_rho_bounded}.
  
  Hence, we can extend the local solutions to global solutions on the whole
  interval $\mathbb{R}_{\ge 0}$ if the condition~\eqref{eqn:L_descent}
  holds almost everywhere.
  In this case,
  equations~\eqref{eqn:norm_y_derivative},~\eqref{eqn:aug_Lagrangian}, and
  Assumption~\ref{ass:phi_rho_bounded} imply that for $t > 0$
  \begin{gather}
    \label{eqn:mean_of_L}
    \begin{aligned}
      \frac{1}{t} \int_{0}^{t} L^{\rho}(x(\tau), y(\tau)) \,\ud \tau
      &= \frac{1}{t} \int_{0}^{t} \phi^{\rho}(x(\tau)) \,\ud \tau
      + \frac{1}{t} \left[ \frac{1}{2} \norm{y(t)}{Y}^2 - \frac{1}{2}
        \norm{y_0}{Y}^2 \right]\\
      &\ge \phi^{\rho}_{\mathrm{low}} - \frac{1}{2t} \norm{y_0}{Y}^2.
    \end{aligned}
  \end{gather}
  Using the monotonicity $L^{\rho}(x(\tau), y(\tau)) \le
  L^{\rho}(x(s), y(s))$ for $0 < s \le \tau$ implied
  by \eqref{eqn:L_descent}, we obtain for $s \le t$ that
  \begin{equation}
    \label{eqn:L_pointwise}
    \frac{1}{t} \int_{0}^{t} L^{\rho}(x(\tau), y(\tau)) \,\ud \tau
    \le 
    \frac{1}{t}
    \int_{0}^{s} L^{\rho}(x(\tau), y(\tau)) \,\ud \tau
    + \frac{t - s}{t} L^{\rho}\left(x(s), y(s)
    \right).
  \end{equation}
  We concatenate \eqref{eqn:mean_of_L} and \eqref{eqn:L_pointwise} and let $t
  \to \infty$, which yields
  \[
    L^{\rho}(x(s), y(s)) \ge \phi^{\rho}_{\mathrm{low}}
    \quad \text{for all } s \in \mathbb{R}_{\ge 0}.
  \]
  Hence, we obtain
  \begin{equation}
    \label{eqn:int_dLdt_bounded}
    0 \ge \int_{0}^{t}\frac{\ud}{\ud \tau} L^{\rho}(x(\tau), y(\tau)) \,\ud\tau
    = L^{\rho}(x(t), y(t)) - L^{\rho}(x_0, y_0)
    \ge \phi^{\rho}_{\mathrm{low}} - L^{\rho}(x_0, y_0).
  \end{equation}
  If condition~\eqref{eqn:c_descent} holds additionally, 
  the boundedness of the integral in \eqref{eqn:int_dLdt_bounded} implies with
  integration of assumption~\eqref{eqn:c_descent} that
  \begin{align}
    &\quad\, 
    \gamma_2 \int_{0}^{t} \norm{c(x(\tau))}{Y}^2 \, \ud \tau
    \label{eqn:int_c_bounded}\\
    &\le -\int_{0}^{t} \frac{\ud}{\ud t} L^{\rho}(x(t), y(t)) \, \ud \tau
    - \gamma_1 \int_{0}^{t} \frac{\ud}{\ud t} \left( \frac{1}{2}
    \norm{c(x(t))}{Y}^2 \right) \ud \tau \nonumber \\
    &\le L^{\rho}(x_0, y_0) - L^{\rho}(x(t), y(t))
    - \gamma_1 \left[ \frac{1}{2} \norm{c(x(t))}{Y}^2 - \frac{1}{2}
    \norm{c(x_0)}{Y}^2 \right] \nonumber \\
    &\le L^{\rho}(x_0, y_0) - \phi^{\rho}_{\mathrm{low}} + \gamma_1
    \frac{1}{2} \norm{c(x_0)}{Y}^2. \nonumber
  \end{align}
  Hence, $\int_{0}^{\infty} \norm{c(x(t))}{Y}^2 \ud t < \infty$ and we can
  establish~\eqref{eqn:bounded_in_L2} by way of~\eqref{eqn:int_dLdt_bounded} and
  the representation~\eqref{eqn:dLdt}.

  If now there is a set $M \subseteq X \times Y$ such that $\nabla
  L^{\rho}$ is Lipschitz continuous on $M$ and $(x(t), y(t)) \in M$ for all $t
  \in [0, \infty)$, then the integrand in~\eqref{eqn:int_c_bounded} is
  absolutely continuous (as a concatenation of an absolutely continuous function
  with Lipschitz continuous functions). This implies uniform continuity of the
  integrand and we can deduce that $\norm{c(x(t))}{Y}^2 \to 0$ for $t \to
  \infty$.  In combination with~\eqref{eqn:int_dLdt_bounded} and the
  representation~\eqref{eqn:dLdt}, this implies that
  \[
    \frac{\ud}{\ud t} L^{\rho}(x(t), y(t))
    = -\norm{P_{T(C, x(t))}\left(-\nabla_x L^{\rho}(x(t), y(t))\right)}{X}^2
    + \norm{c(x(t))}{Y}^2 \to 0
  \]
  and finally $P_{T(C, x(t))}\left( -\nabla_x L^{\rho}(x(t), y(t)) \right) \to
  0$ for $t \to \infty$.
  \qed
\end{proof}

\paragraph{Discussion of Theorem~\ref{thm:existence_of_trajectories}} If we do
not obtain a solution up to $t_{\mathrm{final}} = \infty$, it must be due to
violation of \eqref{eqn:L_descent} or \eqref{eqn:c_descent}. In this case, we
may try to increase $\rho$ in order for the negative term in~\eqref{eqn:dLdt} to
overpower the positive one. To understand the behavior for $\rho \to \infty$, we
let $\beta = 1 / (1+\rho) \in [0, 1]$ and consider a reparametrization of the
flow equations~\eqref{eqn:gradflow} via $x_{\beta}(t) = x(\beta t),
y_{\beta}(t) = y(\beta t)$, which leads to
\begin{align*}
  %\label{eqn:gradflow_beta}
  \dot{x}_{\beta}(t) &= P_{T(C,x_{\beta}(t))} \left(-\beta \nabla_x
    L^{0}(x_{\beta}(t), y_{\beta}(t)) - (1 - \beta) \nabla c(x_{\beta}(t))
  c(x_{\beta}(t)) \right), \\
  \dot{y}_{\beta}(t) &= \beta \nabla_y L^{\rho}(x_{\beta}(t), y_{\beta}(t)).
\end{align*}
For $\beta = 0$, these flow equations reduce to the projected gradient flow for
minimizing the constraint violation $\norm{c(x)}{Y}^2$ over $x \in C$ according
to
\begin{align*}
  %\label{eqn:gradflow_beta}
  \dot{x}_{\beta}(t) &= P_{T(C,x_{\beta}(t))} \left(-\nabla c(x_{\beta}(t))
  c(x_{\beta}(t)) \right), &
  \dot{y}_{\beta}(t) &= 0.
\end{align*}
%For $\beta = 1$, we have that $L^{\rho}(x, y) =
%\frac{1}{2} \norm{c(x(t))}{Y}^2$ and the assumptions \eqref{eqn:L_descent} and
%\eqref{eqn:c_descent} both reduce to
%\begin{align*}
%  \frac{\ud}{\ud t} L^{1}(x(t), y(t)) &= -\norm{P_{T(C, x(t))}\left( -[\nabla
%  c(x(t))] c(x(t)) \right)}{X}^2 \le 0.
%\end{align*}
Hence, violation of \eqref{eqn:L_descent} or \eqref{eqn:c_descent} for large
$\rho$ can only occur if for $\beta = 1$ we get stuck in a
locally infeasible point $\tilde{x}$ of problem \eqref{eqn:mp}, which means
\[
  P_{T(C, \tilde{x})}\left( -[\nabla c(\tilde{x})] c(\tilde{x}) \right) = 0
  \quad \text{but} \quad
  c(\tilde{x}) \neq 0.
\]
This case must arise for instance if $\mathcal{F} = \varnothing$ and it is
reassuring that the theory provides room for this pathological case and that we
at least obtain a point of (locally) minimal constraint violation.

We also remark that boundedness of $y(t)$ can for instance be ensured by the
sufficient condition that for some $\gamma_3 > 0$ we have (omitting
$t$-arguments)
\begin{equation}
  \label{eqn:gronwall_c_descent}
  \begin{aligned}
    \frac{\ud}{\ud t} \left( \frac{1}{2} \norm{c(x)}{Y}^2 \right)
    &= \iprod{P_{T(C, x)} \left(-\nabla_x L^{0}(x, y) - \rho \nabla c(x) c(x)
    \right), \nabla c(x) c(x)}{X}\\
    &\le -\gamma_3 \norm{c(x)}{Y}^2.
  \end{aligned}
\end{equation}
In this case, Gr\"onwall's inequality (see, e.g., \cite{Amann1990}) implies
$\norm{c(x(t))}{Y} \le \norm{c(x_0)}{Y} e^{-\gamma_3 t}$ and consequently
\[
  \norm{y(t) - y_0}{Y} \le \int_{0}^{\infty} \norm{c(x(t))}{Y} \ud t
  \le \gamma_3^{-1} \norm{c(x_0)}{Y}.
\]

Assumption~\eqref{eqn:gronwall_c_descent} is obviously too restrictive for the
case of a feasible initial guess $c(x_0) = 0$, which would imply $y(t) \equiv
y_0$. Hence, we prefer the weaker assumption~\eqref{eqn:c_descent} in
Theorem~\ref{thm:existence_of_trajectories}.

%Another word of caution is appropriate at this point: If the solution exists for
%all $t \in [0, \infty)$, it does not imply that $(x,y)$ converges to a
%stationary point even though $(\dot{x}, \dot{y}) \to 0$, as the following
%example illustrates.
%\begin{example}
%  We consider the unconstrained problem of minimizing $\phi(x) = e^{-x}$ over $x
%  \in \mathbb{R}$, which does not have a minimizer. With the choice $x_0 = 0$,
%  it satisfies Assumption~\ref{ass:phi_rho_bounded} with $\phi_{\mathrm{low}} =
%  0$. Its corresponding gradient flow equation is $\dot{x} = e^{-x}$ with the
%  particular solution $x(t) = \log(1+t)$. Indeed, $\dot{x} = (1+t)^{-1}$ is
%  square integrable in accordance with
%  Theorem~\ref{thm:existence_of_trajectories} but not absolutely integrable on
%  $\mathbb{R}_{\ge 0}$. Thus, the arclength $\int_{0}^{\infty}
%  \abs{\dot{x}(\tau)} \,\ud\tau$ of $x$ is not finite and $x(t) \to \infty$ for
%  $t \to \infty$.
%  In practice, this kind of pathology is seldom an issue.
%\end{example}

We next characterize equilibrium points of \eqref{eqn:gradflow} assuming they
exist.
\begin{lemma}[Equilibria are critical]
  \label{eqn:equilibria_are_critical_points}
  Equilibrium points $(\bar{x}, \bar{y}) \in C \times Y$
  of~\eqref{eqn:gradflow} are critical points of \eqref{eqn:mp} and vice versa.
\end{lemma}
\begin{proof}
  Let $(\bar{x}, \bar{y}) \in C \times Y$ be an equilibrium point of
  \eqref{eqn:gradflow}, implying $0 = \nabla_{y} L^{\rho}(\bar{x},
  \bar{y}) = c(\bar{x})$ and consequently $\bar{x} \in \mathcal{F}.$
  From $0 = P_{T(C, \bar{x})} ( -\nabla_{x} L^{\rho}(\bar{x}, \bar{y}))$, we
  can derive with Lemma~\ref{lem:moreau_decomp} that $- \nabla_{x}
  L^{\rho}(\bar{x}, \bar{y}) = P_{T^{-}(C, \bar{x})}(- \nabla_{x}
  L^{\rho}(\bar{x}, \bar{y})) \in T^{-}(C, \bar{x})$.
  Because $c(\bar{x}) = 0$, we have $\nabla_{x} L^{\rho}(\bar{x}, \bar{y}) =
  \nabla \phi(\bar{x}) + \nabla c(\bar{x}) \bar{y}$. Hence, $(\bar{x},
  \bar{y})$ is a critical point.

  Let now $(\bar{x}, \bar{y}) \in \mathcal{F} \times Y$ be a critical point
  of~\eqref{eqn:mp}. Because $\bar{x} \in \mathcal{F}$, the antigradient
  vanishes due to $\nabla_{y} L^{\rho}(\bar{x}, \bar{y}) = c(\bar{x}) = 0$. By
  definition, we also have that
  \[
    -\nabla_{x} L^{\rho}(\bar{x}, \bar{y})
    = -\nabla \phi(\bar{x}) - \nabla c(\bar{x}) \left[ \bar{y} + \rho c(\bar{x})
    \right]
    = -\nabla \phi(\bar{x}) - \nabla c(\bar{x}) \bar{y} \in T^{-}(C, \bar{x}).
  \]
  Moreau decomposition of $-\nabla_{x} L^{\rho}(\bar{x}, \bar{y})$ then yields
  that $P_{T(C, \bar{x})}(-\nabla_{x} L^{\rho}(\bar{x}, \bar{y})) = 0$. This
  shows that both right-hand sides of~\eqref{eqn:gradflow} vanish and that
  $(\bar{x}, \bar{y})$ is an equilibrium point.
  \qed
\end{proof}

Among the critical points we are apparently only interested in those that are
minima of~\eqref{eqn:mp}. For the finite-dimensional unconstrained case, we
recall that asymptotically stable equilibria of the gradient flow are strict
local minima of the objective function and that the converse is true if the
objective is analytic in a neighborhood of the minimum~\cite{Absil2006}. This is
of high practical relevance, because the gradient flow will be attracted to
strict local minima and, conversely, small perturbations (for instance due to
numerical round-off) will usually make the flow escape unwanted critical points
such as saddle points or maxima.

For the constrained case, the situation is more complicated because the
intrinsic saddle point structure of the Lagrangian requires a
gradient/antigradient flow, for which to our knowledge no results on asymptotic
stability exist so far.
We show that critical points that admit an emanating feasible curve of
descent are not asymptotically stable (under reasonable conditions). This
implies that the projected gradient/antigradient flow will not be attracted to
these undesired critical points. To prove this result, we need the following
three definitions.

\begin{definition}[Descent curve]
  We call a continuous function $\bar{x}: [0, 1] \to \mathcal{F}$ a
  \emph{descent curve} of \eqref{eqn:mp}, if $\phi(\bar{x}(t_2)) <
  \phi(\bar{x}(t_1))$ for all $0 \le t_1 < t_2 \le 1$.
\end{definition}

\begin{definition}[Stability]
  \label{def:stability}
  An equilibrium $(\bar{x}, \bar{y}) \in C \times Y$ of the projected gradient
  flow~\eqref{eqn:gradflow} is \emph{stable} if for every neighborhood $U
  \times V \subset X \times Y$ of $(\bar{x}, \bar{y})$ there exists a smaller
  neighborhood $U_1 \times V_1$ of $(\bar{x}, \bar{y})$ such
  that solutions $(x, y): [0, \infty) \to (U \cap C) \times V$ of
  \eqref{eqn:gradflow} exist for all initial values $(x_0, y_0) \in (U_1 \cap
  C) \times V_1$. If, in addition, it holds for all these solutions that
  $\lim_{t \to \infty} (x(t), y(t)) = (\bar{x}, \bar{y})$, then $(\bar{x},
  \bar{y})$ is \emph{asymptotically stable}.
\end{definition}

\begin{definition}[Flow ribbon]
  \label{def:flow_ribbon}
  For a continuous function $(\bar{x}, \bar{y}): [0, 1] \to C \times Y$ we 
  denote by $\mathcal{R}(\bar{x}, \bar{y}) \subseteq C \times Y$ the \emph{flow
  ribbon emanating from the curve} $(\bar{x}, \bar{y})$, which we define as the
  union of the images of all curves $(x, y): \mathbb{R}_{\ge 0} \to C \times Y$
  satisfying~\eqref{eqn:gradflow} with initial values $(x(0), y(0)) =
  (\bar{x}(l), \bar{y}(l))$ for some $l \in [0, 1]$.
\end{definition}
We can think of a flow ribbon as the trajectory of a curve under the
gradient/antigradient flow~\eqref{eqn:gradflow}, just as if the curve at $t=0$
is the first thread and we weave together the threads into a fabric while moving
along the flow. This somewhat unusual definition is required to keep the set of
points $(x,y)$ small on which assumption~\eqref{eqn:dLdt_lt_0_in_UV} in the
following theorem must hold (compare also Example~\ref{ex:simple_noncvx_qp}
below).

\begin{theorem}
  \label{thm:asymp_stab}
  Let $\bar{x}: [0, 1] \to \mathcal{F}$ be a descent curve and $\bar{y}(t)
  \equiv \bar{y} \in Y$
  such that $(\bar{x}(0), \bar{y})$ is a critical point of \eqref{eqn:mp} and
  let there exist a neighborhood $U \times V \subset X \times Y$ of
  $(\bar{x}(0), \bar{y})$ such that for all $(x, y) \in \mathcal{R}(\bar{x},
  \bar{y}) \cap (U \times V)$ with $L^{\rho}(x, y) < L^{\rho}(\bar{x}(0),
  \bar{y})$ it holds that
  \begin{equation}
    \label{eqn:dLdt_lt_0_in_UV}
    \norm{c(x)}{Y}^2 \le \norm{P_{T(C,x)}\left( -\nabla_x L^{\rho}(x,
    y) \right)}{X}^2.
  \end{equation}
  Then $(\bar{x}(0), \bar{y})$ is not asymptotically stable.
\end{theorem}
\begin{proof}
  \emph{\negmedspace by contradiction.}
  Assume $(\bar{x}(0), \bar{y})$ is asymptotically stable.
  By Definition~\ref{def:stability}, there exists a neighborhood $U_1 \times V_1
  \subset U \times V$ of $(\bar{x}(0), \bar{y})$, which admits for each element
  as initial value a global solution to \eqref{eqn:gradflow}.  We choose $l \in
  (0, 1]$ such that $(x_0, y_0) := (\bar{x}(l), \bar{y}) \in U_1 \times V_1$.
  Because $\bar{x}$ is a descent curve, we have that $c(x_0) = 0$ and
  \begin{equation}
    \label{eqn:L0_lt_Lbar}
    L^{\rho}(\bar{x}(0), \bar{y}) - L^{\rho}(x_0, y_0)
    = \phi(\bar{x}(0)) - \phi(x_0) =: \varepsilon > 0.
  \end{equation}
  By Definition~\ref{def:stability}, a solution $(x, y): [0, \infty) \to (U \cap
  C) \times V$ of \eqref{eqn:gradflow} with $x(0) = x_0$ and $y(0) = y_0$
  exists and converges to $(\bar{x}(0), \bar{y})$. Using
  assumption~\eqref{eqn:dLdt_lt_0_in_UV} and equations~\eqref{eqn:dLdt}
  and~\eqref{eqn:L0_lt_Lbar}, we observe that
  \begin{align*}
    L^{\rho}(x(t), y(t)) &= L^{\rho}(x_0, y_0) + \int_{0}^{t}
    \frac{\ud}{\ud t} L^{\rho}(x(\tau), y(\tau)) \,\ud\tau\\
    &\le L^{\rho}(x_0, y_0) = L^{\rho}(\bar{x}(0), \bar{y}) - \varepsilon
  \end{align*}
  for all $t \in \mathbb{R}_{\ge 0}$,
  which implies that $(x, y)$ cannot converge to $(\bar{x}(0), \bar{y})$. Hence,
  $(\bar{x}(0), \bar{y})$ is not asymptotically stable.
  \qed
\end{proof}

In order to validate that assumption~\eqref{eqn:dLdt_lt_0_in_UV} does not reduce
the assertion of Theorem~\ref{thm:asymp_stab} to one about the empty set, we
provide a simple example.

\begin{example}[Simple nonconvex quadratic program]
  \label{ex:simple_noncvx_qp}
  We consider the problem
  \[
    \min \tfrac{1}{2} \left(x_1^2 -x_2^2\right) \quad
    \text{over } x \in \mathbb{R} \times \mathbb{R}_{\ge 0} =: C \quad
    \text{subject to } x_1 = 0.
  \]
  It is easy to verify that $(x_1, x_2, y) = (0, 0, 0)$ is a critical point and
  the objective is unbounded for the feasible points $x_1 = 0$, $x_2 \to
  \infty$. The augmented Lagrangian amounts to
  \[
    L^{\rho}(x,y) = \tfrac{1}{2} (x_1^2 - x_2^2) + x_1 y + \tfrac{\rho}{2} x_1^2
    = (1+\rho)\tfrac{1}{2} x_1^2 - \tfrac{1}{2} x_2^2 + x_1 y.
  \]
  The projected gradient flow equations then read (omitting $(t)$-arguments)
  \begin{align*}
    \dot{x} &= P_{T(C,x)}(-\nabla_x L^{\rho}(x,y)) =
    \begin{pmatrix}
      -(1+\rho) x_1 - y\\
      \max(0, x_2)
    \end{pmatrix},
    &
    \dot{y} &= \nabla_y L^{\rho}(x, y) = x_1.
  \end{align*}
  For the descent curve $\bar{x}(l) = (0, l, 0)^T, l \in [0, 1],$ with
  corresponding $\bar{y} \equiv 0$, we can easily solve the flow equations and
  obtain the flow ribbon
  \[
    \mathcal{R}(\bar{x},\bar{y}) = \left\{
      (0, l e^t, 0)^T \mid t \in \mathbb{R}, l \in [0, 1] \right\}
    = \{0\} \times \mathbb{R}_{\ge 0} \times \{0\}.
  \]
  Hence, we see that for $(x, y) \in \mathcal{R}(\bar{x}, \bar{y})$ it holds
  that
  \[
    L^{\rho}(x, y) = -\tfrac{1}{2} x_2^2 \le 0 = L^{\rho}(\bar{x}(0), \bar{y})
  \]
  with strict inequality for $x_2 > 0$. Inequality~\eqref{eqn:dLdt_lt_0_in_UV}
  holds for all $(x, y) \in \mathcal{R}(\bar{x}, \bar{y})$ by virtue of
  \[
    \norm{c(x)}{2}^2 = x_1^2 = 0 \le \max(0, x_2)^2
    = \norm{P_{T(C, x)}(-\nabla_x L^{\rho}(x, y))}{2}^2.
  \]
  Hence, this example satisfies all assumptions of Theorem~\ref{thm:asymp_stab}.
\end{example}

\section{Projected backward Euler: A sequential homotopy method}
\label{sec:seqhom}

It is well-known that the projection in \eqref{eqn:gradflow} is
actually the derivative of the projection of the primal variable onto $C$ in
direction of the negative primal gradient:
\begin{lemma}
  \label{lem:projection_gateaux_diff}
  For a nonempty closed convex set $K \subseteq X$,
  the G\^{a}teaux derivative of the projection of $x \in X$ onto $K$ in the
  direction $\delta x \in X$ is the projection of $\delta x$ onto the tangent
  cone $T(K, x)$, i.e.,
  \[
    %\Pi_{K}(x, \delta x) :=
    \lim_{h \to 0^+} h^{-1} \left( P_{K}(x + h \delta x) - x \right)
    = P_{T(K,x)}(\delta x).
  \]
\end{lemma}
\begin{proof}
  See \cite[Lemma~4.5]{Zarantonello1971}. \qed
\end{proof}
This motivates following the flow defined by~\eqref{eqn:gradflow} from
$(\hat{x}, \hat{y}) \in C \times Y$ to $(x, y) \in C \times Y$ with a projected 
backward Euler step of step size $\Delta t > 0$ by solving
\begin{align}
  \label{eqn:backwardEuler}
  x - P_{C}\left( \hat{x} - \Delta t \nabla_x L^{\rho}(x, y) \right) &= 0, &
  y - \hat{y} - \Delta t c(x) &= 0,
\end{align}
because Lemma~\ref{lem:projection_gateaux_diff} ensures consistency by virtue of
\[
  \lim_{\Delta t \to 0} \frac{x - \hat{x}}{\Delta t}
  = \lim_{\Delta t \to 0} \frac{P_{C}\left( \hat{x} - \Delta t \nabla_x
  L^{\rho}(x, y) \right) - \hat{x}}{\Delta t}
  = P_{T(C, \hat{x})} \left( -\nabla_x L^{\rho}(\hat{x}, \hat{y}) \right).
\]
From a computational point of view, the projected backward Euler system
\eqref{eqn:backwardEuler} is an ideal candidate for the application of local
(possibly inexact) semismooth Newton methods (see, e.g.,
\cite{Mifflin1977,Qi1993,Ulbrich2011}), which we will investigate in more detail
in section~\ref{sec:PDE_opt}.

In addition, the projected backward Euler system \eqref{eqn:backwardEuler} can
be interpreted as necessary optimality conditions of a primal-dual proximally
regularized version of the augmented form of~\eqref{eqn:mp}. With $\lambda =
1/\Delta t$, it reads
\begin{equation}
  \label{eqn:homotopy_mp}
  \begin{alignedat}{3}
    \min~ & \phi^{\rho}(x) + \lambda \left[ \tfrac{1}{2} \norm{x -
    \hat{x}}{X}^2 + \tfrac{1}{2} \norm{w - \hat{y}}{Y}^2 \right]
    & ~~\text{over}~&w \in Y, x \in C\\
    \text{subject to}\,~ & c(x) + \lambda w = 0.
  \end{alignedat}
\end{equation}
Uniqueness of solutions to~\eqref{eqn:homotopy_mp} can be guaranteed for
sufficiently large $\lambda$.
\begin{theorem}
  \label{thm:homotopy_mp}
  The regularized problem~\eqref{eqn:homotopy_mp} has the following
  properties for $\lambda > 0$:
  \begin{enumerate}
    \item It satisfies the strong constraint qualification of
      Lemma~\ref{lem:IFCQ}.
    \item Its primal-dual solutions $(\bar{w}, \bar{x}, \bar{y}) \in Y \times C
      \times Y$ satisfy \eqref{eqn:backwardEuler} and $\bar{w} = -\Delta t
      c(\bar{x})$.
    \item For $\lambda \to \infty$, i.e., $\Delta t \to 0$, its unique
      primal-dual solution $(\bar{w}, \bar{x}, \bar{y})$ tends to $(0, \hat{x},
      \hat{y})$ provided that $\nabla L^{\rho}$ is globally Lipschitz
      continuous.
  \end{enumerate}
  For $\lambda = 0$, i.e., $\Delta t = \infty$, its primal-dual solutions
  $\bar{x}$ and $\bar{y}$ coincide with those of problem~\eqref{eqn:mp} and
  arbitrary $\bar{w} \in Y$.
\end{theorem}
\begin{proof}
  With $U = Y$, $Q = X$, we can apply Lemma~\ref{lem:IFCQ} with the solution
  mapping $S(x) = -\Delta t c(x)$ and $\range \lambda \eye_{Y} = \range \lambda
  \eye_{Y}^{\star} = Y$. This shows assertion 1.
  We call the Lagrangian of
  \eqref{eqn:homotopy_mp} homotopy Lagrangian or proximal Lagrangian and denote
  it by
  \[
    L^{\lambda,\rho}(w, x, y) = L^{\rho}(x, y) + \lambda \left[
    \tfrac{1}{2} \norm{x - \hat{x}}{X}^2 + \tfrac{1}{2} \norm{w -
    \hat{y}}{Y}^2 + \iprod{y, w}{Y} \right].
  \]
  By Lemma~\ref{lem:IFCQ}, GCQ holds at all feasible points and
  Theorem~\ref{thm:noc} yields that
  \begin{equation}
    \label{eqn:homotopy_stationary}
    \left(-\nabla_w L^{\lambda,\rho}(\bar{w},\bar{x},\bar{y}),
    -\nabla_x L^{\lambda,\rho}(\bar{w},\bar{x},\bar{y}) \right)
    \in \{0\} \times T^{-}(C, \bar{x}),
  \end{equation}
  from which we can deduce that $\bar{w} = \hat{y} - \bar{y}$ because of
  \begin{equation}
    \label{eqn:stat_w}
    0 = \nabla_w L^{\lambda, \rho} (\bar{w}, \bar{x}, \bar{y}) = 
    \lambda (\bar{w} - \hat{y} + \bar{y}).
  \end{equation}
  Hence, the feasibility of $(\bar{w}, \bar{x})$ implies that
  \[
    c(\bar{x}) + \lambda [\hat{y} - \bar{y}] = 0.
  \]
  Multiplication with $\Delta t$ yields the second equation
  of~\eqref{eqn:backwardEuler}. For the $x$-part of
  \eqref{eqn:homotopy_stationary}, we observe that
  \[
    -\Delta t \nabla_x L^{\lambda,\rho} (\bar{w}, \bar{x}, \bar{y}) =
    -\bar{x} + \hat{x} - \Delta t \nabla_x L^{\rho}(\bar{x}, \bar{y}) \in
    T^{-}(C, \bar{x}),
  \]
  implying $\hat{x} - \Delta t \nabla_x L^{\rho}(\bar{x}, \bar{y}) \in T^{-}(C,
  \bar{x}) + \bar{x}$
  and by Lemma~\ref{lem:polar_cone_and_projection} that therefore $P_{C}(\hat{x}
  - \Delta t \nabla_x L^{\rho}(\bar{x}, \bar{y})) = \bar{x}$, which coincides
  with the first equation of~\eqref{eqn:backwardEuler}. This shows assertion 2.

  We can now use~\eqref{eqn:backwardEuler} to define a fixed point iteration
  $z^{k+1} = \Phi(z^k)$ on $C \times Y$ via
  \[
    \Phi(z) = (P_C(\hat{x} - \Delta t \nabla_{x} L^{\rho}(z)),
    \hat{y} + \Delta t \nabla_{y} L^{\rho}(z)).
  \]
  Let $\omega$ denote the Lipschitz constant of $\nabla L^{\rho}$. For $\Delta t
  < \frac{1}{\omega}$, the mapping $\Phi$ is a contraction because $P_C$ is
  Lipschitz continuous with modulus 1: 
  \begin{multline*}
    \norm{\Phi(z) - \Phi(\tilde{z})}{X \times Y}^2
    = \norm{P_C(\hat{x} - \Delta t \nabla_{x} L^{\rho}(z)) - P_C(\hat{x} -
    \Delta t \nabla_{x} L^{\rho}(\tilde{z}))}{X}^2\\
    + \Delta t^2 \norm{\nabla_{y} L^{\rho}(z) - \nabla_{y}
    L^{\rho}(\tilde{z})}{Y}^2
    \le (\omega \Delta t)^2 \norm{z - \tilde{z}}{X \times Y}^2.
  \end{multline*}
  The Banach fixed point theorem yields uniqueness
  and existence of a fixed point $(\bar{x}, \bar{y})$, which together with $\bar{w}
  = \hat{y} - \bar{y}$ is the unique solution of~\eqref{eqn:homotopy_mp}. 
  For $\Delta t = 0$, the fixed point and thus the solution is obviously $(0,
  \hat{x}, \hat{y})$.
  In order to prove convergence of $(\bar{w}, \bar{x}, \bar{y})$ to $(0,
  \hat{x}, \hat{y})$ for $\Delta t \to 0$, we observe that
  \[
    \norm{\bar{z} - \hat{z}}{X \times Y} = 
    \norm{\Phi(\bar{z}) - \Phi(\hat{z}) + \Phi(\hat{z}) - \hat{z}}{X \times Y}
    \le \omega \Delta t \norm{\bar{z} - \hat{z}}{X \times Y}
    + \norm{\Phi(\hat{z}) - \hat{z}}{X \times Y},
  \]
  which implies (recall that $\Phi(\hat{z})$ depends continuously on $\Delta t$)
  \[
    \norm{\bar{z} - \hat{z}}{X \times Y}
    \le \frac{1}{1 - \omega \Delta t} \norm{\Phi(\hat{z}) - \hat{z}}{X \times Y}
    \to 0 \quad \text{for } \Delta t \to 0.
  \]
  This finally proves assertion 3.
  \qed
\end{proof}

The artificial introduction of the variable $w$ in \eqref{eqn:homotopy_mp}
allows a lifting of the dual regularization term $y - \hat{y}$ in the backward
Euler system \eqref{eqn:backwardEuler} onto primal variables. From a linear
algebra perspective, this can be understood as a Schur complement approach, as
we see in the following example.

\begin{example}[A quadratic program]
  Let $C = X$, $\rho = 0$, $\Delta t > 0$, $c(x) = A x - b$, and $\phi(x) =
  \frac{1}{2} \iprod{x, H x}{X} - \iprod{g, x}{X}$ for $A \in \mathcal{L}(X,
  Y)$, $H = H^{\star} \in \mathcal{L}(X, X)$, $g \in X$, and $b \in Y$. The
  necessary optimality conditions of the homotopy
  problem~\eqref{eqn:homotopy_mp} are then equivalent to the linear system
  \[
    \begin{pmatrix}
      \lambda \eye_{Y} & 0 & \lambda \eye_{Y} \\
      0 & H + \lambda \eye_{X} & A^{\star}\\
      \lambda \eye_{Y} & A & 0
    \end{pmatrix}
    \begin{pmatrix}
      \bar{w}\\
      \bar{x}\\
      \bar{y}
    \end{pmatrix}
    =
    \begin{pmatrix}
      \lambda \hat{y}\\
      \lambda \hat{x} + g\\
      b
    \end{pmatrix}.
  \]
  If we eliminate $\bar{w}$ with a Schur complement approach, we obtain the
  backward Euler system~\eqref{eqn:backwardEuler} as a primal-dual
  regularization of the original saddle point system for \eqref{eqn:mp}
  according to
  \[
    \left[
      \begin{pmatrix}
        H & A^{\star}\\
        A & 0
      \end{pmatrix}
      + \lambda
      \begin{pmatrix}
        \eye_{X} & 0\\
        0 & -\eye_{Y}
      \end{pmatrix}
    \right]
    \begin{pmatrix}
      \bar{x}\\
      \bar{y}
    \end{pmatrix}
    =
    \begin{pmatrix}
      g + \lambda \hat{x}\\
      b - \lambda \hat{y}
    \end{pmatrix}.
  \]
\end{example}

We can derive two interesting equivalent reformulations
of~\eqref{eqn:homotopy_mp}. The first reformulation substitutes
$v = \sqrt{\lambda} w$, from which we obtain
\begin{equation}
  \label{eqn:homotopy_mp_sqrt}
  \begin{alignedat}{3}
    \min~ & \phi^{\rho}(x) + \frac{\lambda}{2} \norm{x - \hat{x}}{X}^2
    + \frac{1}{2} \norm{v - \sqrt{\lambda} \hat{y}}{Y}^2
    & ~~\text{over}~&v \in Y, x \in C\\
    \text{subject to}\,~ & c(x) + \sqrt{\lambda} v = 0.
  \end{alignedat}
\end{equation}
The advantage of~\eqref{eqn:homotopy_mp_sqrt} over~\eqref{eqn:homotopy_mp} is
that the optimal $v$ is also uniquely determined for $\lambda = 0$.
The second reformulation completely eliminates $w = -\Delta t c(x)$. This leads
to the problem
\[
  \begin{alignedat}{3}
    \min~ & \phi^{\rho}(x) + \frac{\lambda}{2} \norm{x - \hat{x}}{X}^2 +
    \frac{1}{2} \norm{\sqrt{\lambda} \hat{y} + \sqrt{1/\lambda} c(x)}{Y}^2
    & ~~\text{over}~& x \in C,
  \end{alignedat}
\]
which has no equality constraint and might allow for the application of
projected Newton/gradient methods similar to,
e.g.,~\cite{Bertsekas1982,Calamai1987,Kelley1994}.

The homotopy problem~\eqref{eqn:homotopy_mp} and Theorem~\ref{thm:homotopy_mp}
provide a complementary interpretation of using projected backward Euler
steps~\eqref{eqn:backwardEuler} for the gradient flow
equations~\eqref{eqn:gradflow}: We trace the solutions
of~\eqref{eqn:homotopy_mp} from some primal-dual starting point $(0, \hat{x},
\hat{y})$ as a continuation in $\lambda$ until the homotopy breaks down. The
result yields an update for $(\hat{x}, \hat{y})$ and we can repeat the
procedure. If, at one point, we are able to drive $\lambda$ to
zero, we can solve the original problem~\eqref{eqn:mp} with superlinear
local convergence rate by the means of a locally superlinearly convergent method
for the homotopy problem~\eqref{eqn:homotopy_mp}, e.g., a semismooth Newton
method. If it is
never possible to drive $\lambda$ to zero, we at least follow the gradient
flow~\eqref{eqn:gradflow} with a projected backward Euler method with stepsize
$1/\lambda$.
If we fix $\lambda$ to some positive value, we obtain a locally linear
convergence rate provided that the gradient flow converges exponentially.

\section{Numerical case study in PDE constrained optimization}
\label{sec:PDE_opt}

We apply the proposed method to the following benchmark problem adapted from
\cite{Lubkoll2017}: Let $\Omega \subset \mathbb{R}^2$ be a bounded domain with
Lipschitz boundary and let constants $a, b, \gamma > 0$, control bounds
$q_{\mathrm{l}}, q_{\mathrm{u}} \in L^{r}(\Omega), r \in (2, \infty]$, and a
target function $u_{\mathrm{d}} \in L^2(\Omega)$ be given. We solve the
control-constrained quasilinear elliptic optimal control problem
\begin{equation}
  \label{eqn:model_problem}
  \begin{alignedat}{3}
    \min~&\frac{1}{2} \int_{\Omega} \abs{u - u_{\mathrm{d}}}^2 %\ud\xi
    + \frac{\gamma}{2} \int_{\Omega} \abs{q}^2 \quad %\ud\xi
    && \text{over } u \in H^1_0(\Omega), q \in L^2(\Omega)\\
    \text{subject to} ~& \nabla \cdot \left(\left[a + b \abs{u}^2\right] \nabla
    u\right) = q,\\
    & q_{\mathrm{l}} \le q \le q_{\mathrm{u}}.
  \end{alignedat}
\end{equation}
In addition to \cite{Lubkoll2017}, we include pointwise control bounds. For
smaller values of $a$ and $\gamma$, problem~\eqref{eqn:model_problem} becomes
more and more ill-conditioned, while the effects of nonlinearity become more
challenging for larger values of $b$.

To transform problem~\eqref{eqn:model_problem} into the form~\eqref{eqn:mp}, we
use the variables $x = (u, q) \in X = U \times Q = H^1_0(\Omega) \times
L^2(\Omega)$, $y \in Y = U^{\ast} = H^{-1}(\Omega)$ and define the closed convex
set
\[
  C = U \times
  \left\{ q \in Q \mid q_{\mathrm{l}} \le q \le q_{\mathrm{u}} \right\}
  =: U \times C_{Q} \subset U \times Q = X
\]
and the functions $\phi: X \to \mathbb{R}$ and $c: X \to Y$ via
\begin{align*}
  \phi((u,q)) &= \frac{1}{2} \int_{\Omega} \abs{u - u_{\mathrm{d}}}^2 %\ud\xi
  + \frac{\gamma}{2} \int_{\Omega} \abs{q}^2\\ %\ud\xi\\
  \dpair{c((u,q)), \varphi}{U} &= \int_{\Omega} \nabla
  \varphi \cdot \left[ a + b \abs{u}^2 \right] \nabla u %\,\ud\xi
  - \int_{\Omega} \varphi q %\,\ud\xi
  \quad \text{for all } \varphi \in U,
\end{align*}
where $c$ is the weak form of the PDE in~\eqref{eqn:model_problem}.
The problem has a continuously Fr\'{e}chet-differentiable solution operator $S:
Q \to U$ in the sense of Lemma~\ref{lem:IFCQ}~\cite{Casas2009}.

\subsection{Implementation aspects}

From an implementation point of view, the projected backward Euler
system~\eqref{eqn:backwardEuler} with all its required derivatives can be
conveniently generated by the use of the Unified Form
Language~\cite{AlnaesEtAl2012,Alnaes2012a} in combination with
Algorithmic Differentiation~\cite{Griewank2008}, as it is implemented in the
DOLFIN/FEniCS
project~\cite{LoggWells2010a,LoggWellsEtAl2012a,AlnaesBlechta2015a,LoggMardalEtAl2012a}.

When evaluating the augmented objective $\phi^{\rho}(x) = \phi(x) +
\frac{\rho}{2} \norm{c(x)}{Y}^2$, the inner product $\iprod{y, c(x)}{Y}$, or the
dual proximal term in~\eqref{eqn:homotopy_mp}, we face the problem of computing
norms and inner products in $Y=H^{-1}(\Omega)$, which we can facilitate
computationally with the
use of the Riesz isomorphism $\norm{y}{Y} = \norm{R_{U} y}{U}$. If we choose the
norm $\norm{u}{U} = \norm{\nabla u}{L^2(\Omega)^2}$ on $U$, the evaluation of
$R_{U} y$ boils down to one solution of a Poisson problem with right-hand side
$y$ and homogeneous Dirichlet boundary conditions. The difficulty from a
computational vantage point is that $R_{U}$ is a large dense matrix in contrast
to its inverse $R_{U}^{-1}$, which is a sparse finite element stiffness matrix.
For practical purposes, we always work with the Riesz represenation of the dual
variable $y_{R} = R_{U} y$ directly, eliminating the need for evaluating the
Riesz isomorphism for the dual variables. 

From a linear algebra point of view, it is important to exploit the special
structure of the augmentation term $\frac{\rho}{2} \norm{c(x)}{Y}^2$. We extend
a well-known argument for the special case of $\lambda = 0$ (see,
e.g.,~\cite[p.~158f]{Ito2008}) to the case $\lambda \ge 0$:
For fixed $(x, y)$, let us denote the gradients and the second derivative of the
augmented Lagrangian $L^{\rho}(x, y)$ by
\begin{align*}
  \nabla_{x} L^{\rho}(x, y) &= \nabla_{x} L^{0}(x, y + \rho c(x)) =: F_1, \quad
  \nabla_{y} L^{\rho}(x, y) = c(x) =: F_2,\\
  \nabla_{xx} L^{\rho}(x, y) &= \nabla_{xx} L^{0}(x, y + \rho c(x)) + \rho
  \nabla c(x) c'(x) =: H + \rho A^{\star} A.
\end{align*}
Disregarding inequalities for a moment, each Newton step for the (appropriately
scaled) backward Euler equations~\eqref{eqn:backwardEuler} requires us to solve
the linear system
\begin{equation}
  \label{eqn:augmented_Newton_system}
  \begin{pmatrix}
    \lambda \eye_{X} + H + \rho A^{\star} A & A^{\star}\\
    A & -\lambda \eye_{Y}
  \end{pmatrix}
  \begin{pmatrix}
    \delta x\\ \delta y
  \end{pmatrix}
  = -
  \begin{pmatrix}
    F_1 + \lambda (x - \hat{x})\\
    F_2 - \lambda (y - \hat{y})
  \end{pmatrix}.
\end{equation}
The problem here is that $A^{\star}A = R_{X} A^{\ast} R_{Y}^{-1} A =
R_{X} A^{\ast} R_{U} A$ becomes a dense matrix after discretization by finite
elements due to $R_{U}$. Hence, we must avoid the formation of $A^{\star} A$.
Instead of~\eqref{eqn:augmented_Newton_system} we solve the equivalent system
\begin{equation}
  \label{eqn:Newton_system}
  \begin{pmatrix}
    \lambda \eye_{X} + H & A^{\star}\\
    A & -(1 + \rho \lambda)^{-1} \lambda \eye_{Y}
  \end{pmatrix}
  \begin{pmatrix}
    \delta x\\
    \delta \tilde{y}
  \end{pmatrix}
  = -
  \begin{pmatrix}
    F_1 + \lambda (x - \hat{x})\\
    (1 + \rho \lambda)^{-1} \left( F_2 - \lambda (y - \hat{y}) \right)
  \end{pmatrix}
\end{equation}
with the reconstruction $\delta y = (1 + \rho \lambda)^{-1} (\delta \tilde{y} +
\rho F_2)$. The equivalence can easily be checked. Because we work with
$y_{R} = R_{U} y$ directly, we need to compute the Riesz representation
$c_{R} = R_{U} c(x)$ first, evaluate the Lagrangian derivatives at $(x,
y_{R} + \rho c_{R})$, solve the unaugmented Newton
system~\eqref{eqn:Newton_system} (reformulated for $y_{R}$ instead of $y$) for
$(\delta x, \delta \tilde{y}_{R})$, and finally reconstruct $\delta y_{R} = (1 +
\rho \lambda)^{-1} ( \delta \tilde{y}_{R} + \rho c_{R})$.

The enforcement of the projection onto $C$ in~\eqref{eqn:backwardEuler}
can be easily implemented on top of~\eqref{eqn:Newton_system}: Let us consider
the block row corresponding to the gradient with respect to $u$
in~\eqref{eqn:backwardEuler} scaled by $\lambda$, which reads
\begin{align*}
  0 &=
  \lambda q - \lambda P_{C_{Q}} \left( \hat{q} - \Delta t \nabla_{q}
  L^{\rho}( (u,q), y) \right)\\
  &= \lambda q - \lambda P_{C_{Q}} \left( \hat{q} - \Delta t \left[ \gamma q -
  R_{U} (y + \rho c(x)) \right] \right).
\end{align*}
This nonsmooth equation together with the remaining smooth block rows
of~\eqref{eqn:backwardEuler} scaled by $\lambda$ can be solved
efficiently with a semismooth Newton method. To this end, we need to address a
norm gap for the pointwise defined projector
\[
  P_{C_{Q}}(q)(\xi) = \max(q_{\mathrm{l}}(\xi), \min(q(\xi),
  q_{\mathrm{u}}(\xi))) \quad \text{for } \xi \in \Omega,
\]
which is known to be semismooth only if it maps from $L^{r}(\Omega) \subsetneq
Q$ to $Q = L^2(\Omega)$ (see, e.g., \cite[sec.~3.3]{Ulbrich2011} or
\cite[Theorem~4.2]{Hintermueller2010}). Indeed, this higher regularity holds
here if the initial guess satisfies $q_0 \in L^{r}(\Omega)$: For
problem~\eqref{eqn:model_problem}, the $Q$ part of the projected backward Euler
equations~\eqref{eqn:backwardEuler} simplify to
$q = P_{C_{Q}} \left( \hat{q} - \Delta t \left[ \gamma q - R_{U} (y + \rho
  c(x)) \right] \right).$
By induction, we can assume that $q, \hat{q} \in L^{r}(\Omega)$. Then, the
argument of the projection operator $P_{C_{Q}}$ also lies in $L^{r}(\Omega)$,
because $R_{U} [y + \rho c(x)] \in H^{1}_{0}(\Omega)$, which is continuously
embedded in $L^{r}(\Omega)$. Because $q_{\mathrm{l}}, q_{\mathrm{u}} \in
L^{r}(\Omega)$, we obtain $q = P_{C_{Q}}(.) \in L^{r}(\Omega)$, which completes
the induction step.

\subsection{Solution algorithm}

\begin{algorithm}[tb]
  \caption{Sequential homotopy method}
  \label{alg:seqhom}
  \DontPrintSemicolon
  \KwData{$z_0 = (x_0, y_0) \in C \times Y = Z$, $\Theta \in (0, 1)$,
    $\lambda_{\mathrm{term}} > 0$, $\lambda_{\mathrm{inc}} > 1$,
    $\mathrm{TOL} > 0$}
  Initialize $z = (x, y) \leftarrow z_0$, $\lambda \leftarrow 1$\;
  \Loop{(over homotopies)}{
    Initialize $\hat{z} = (\hat{x}, \hat{y}) \leftarrow z$\;
    \Loop{(to trace single homotopy leg)}{
      Compute $z^{+}$ by one semismooth Newton step for
      \eqref{eqn:backwardEuler} starting from $z$\;
      \label{ln:Newton_step}
      Compute $z^{++}$ by one simplified semismooth Newton step
      for \eqref{eqn:backwardEuler} starting from $z^{+}$\;
      \label{ln:simplified_Newton_step}
      \uIf{$\norm{z^{++} - z^{+}}{Z} \le \Theta \norm{z^{+} - z}{Z}$}{
        \label{ln:monotonicity_test}
        Accept iterate $z \leftarrow z^{++}$\;
        \lIf{$\lambda \le \lambda_{\mathrm{term}}$ and $\norm{z - \hat{z}}{Z}
          \le \mathrm{TOL}$}{
          \Return solution $z$
        }
        Update homotopy parameter $\lambda$ (e.g., by a proportional-integral
        controller)\;
        \Break inner loop
      }
      \lElse{
        increase homotopy parameter $\lambda \leftarrow \lambda_{\mathrm{inc}}
        \lambda$
      }
    }
  }
\end{algorithm}

We provide in Algorithm~\ref{alg:seqhom} pseudocode for a prototypical
implementation of the sequential homotopy method with a classical continuation
approach. It consists of an outer loop over the subsequent homotopies. In the
inner loop, the reference point $\hat{z} = (\hat{x}, \hat{y})$ is fixed and we
trace the solution of~\eqref{eqn:backwardEuler} with one semismooth Newton step
followed by one inexact semismooth Newton step.

The computationally heavy part is the computation of $z^{+}$ in
line~\ref{ln:Newton_step} by one local semismooth Newton step at $z$ and of
$z^{++}$ in line~\ref{ln:simplified_Newton_step} by one local simplified
semismooth Newton step at $z^{+}$. Here, \emph{simplified} means that the
system matrix of the previous semismooth Newton system is reused, subject to
modifications concerning the current active set guess derived from the residual
evaluated at $z^{+}$.  We accept an iterate for the current value of $\lambda$
if the following natural monotonicity test is satisfied in
line~\ref{ln:monotonicity_test}: We require that the simplified semismooth
Newton increment is smaller in norm than a contraction factor $\Theta \in (0,
1)$ times the semismooth Newton increment.

If the monotonicity test fails, we enlarge $\lambda$ by a constant factor to
drive the solution of~\eqref{eqn:backwardEuler} closer to $\hat{z}$ in order to
eventually enter the region of local superlinear convergence of the semismooth
Newton method.

If the monotonicity test is satisfied, we accept $z^{++}$ as the new iterate.
If $\lambda$ and the norm of the outer loop increment $z - \hat{z}$ are small
enough, then we terminate with the solution $z$, otherwise we predict a new step
size which should eventually drive $\lambda$ close to zero. 
We then commence the next outer iteration.

There are many possibilities to predict the next $\lambda$ after acceptance of
the current iterate. For the numerical results below, we use a heuristic
motivated by a discrete proportional-integral (PI) controller: We try to choose
$\lambda$ such that the contraction factor $\theta = \norm{z^{++} - z^{+}}{Z} /
\norm{z^{+} - z}{Z}$ is close to a given reference $\theta_{\mathrm{ref}} \in
(0, 1)$. We choose to predict $\lambda \leftarrow \lambda /
\lambda_{\mathrm{mod}}$, where $\log \lambda_{\mathrm{mod}}$ is the manipulated
variable. To this end, let $e = \log \theta_{\mathrm{ref}} - \log \theta$ and
let $I$ denote the sum of all previous values of $e$ over the last successful
outer loops. We then set with some constants $K_{P}$ and $K_{I}$
\[
  \log \lambda_{\mathrm{mod}} \leftarrow K_{P} e + K_{I} I.
\]
In each accepted iteration, we have the simple update $I \leftarrow I + e$.
In case the monotonicity test fails, we possibly reset the integral term $I
\leftarrow \min(I, 0)$. We can also clip $\lambda$ at a lower bound
$\lambda_{\mathrm{min}}$. For a
related concept in the stepsize control of one-step methods for ordinary
differential equations we refer to~\cite[p.~28ff]{Hairer1996}.

It is also possible to keep all iterates inside $C$ with an additional
projection in the local semismooth Newton step (see, e.g.,~\cite{Ulbrich2011}).
We found the method to require fewer iterations on~\eqref{eqn:model_problem}
without projection steps, even though we are aware that if $z \not \in C$, we
might run into problems with the monotonicity test in
line~\ref{ln:monotonicity_test} of Algorithm~\ref{alg:seqhom} because
$\norm{z^{+} - z}{Z}$ might not tend to 0 for $\lambda \to \infty$.

Alternatively to Algorithm~\ref{alg:seqhom}, it is conceivable to update the
reference point $\hat{z}$ less frequently and to trace each homotopy leg until
it nearly breaks down in a singularity. In our experience, this approach of long
homotopy legs leads to a more complicated algorithm and requires the solution of
more and worse conditioned linear systems. We prefer the sequential homotopy
method with short homotopy legs in the form of Algorithm~\ref{alg:seqhom}.

\subsection{Numerical results}

\begin{figure}[tb]
  \centering
  \includegraphics[width=\textwidth]{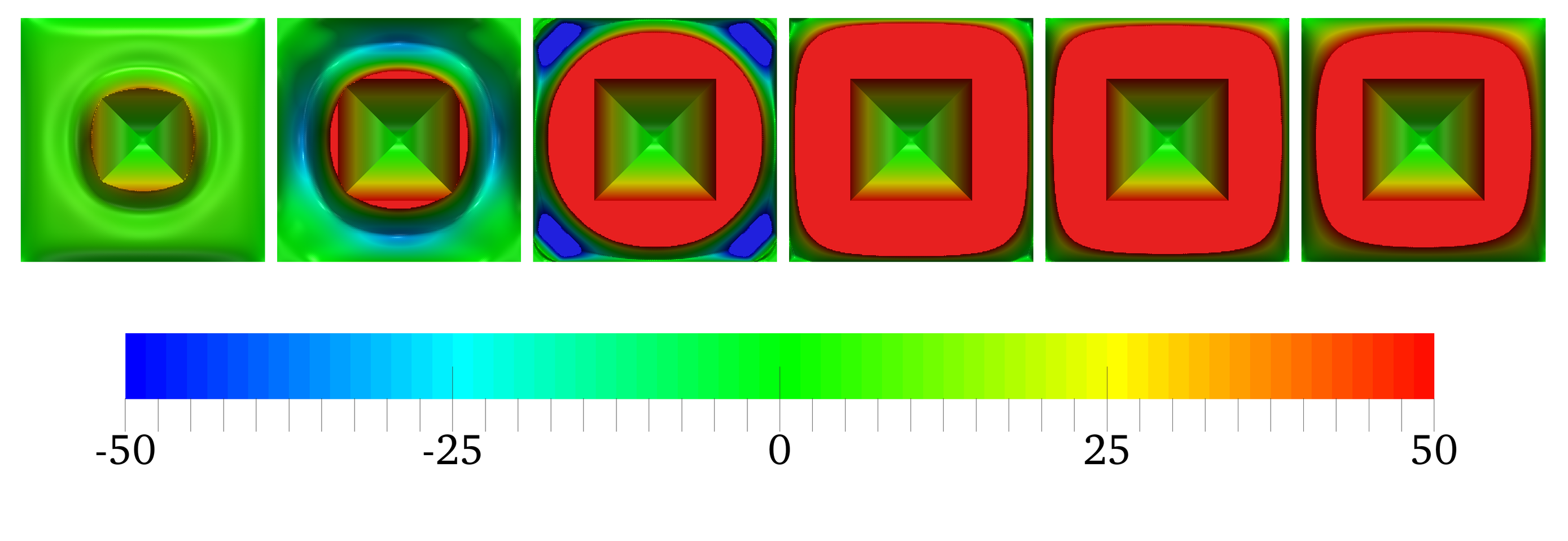}
  \caption{Optimal controls for problem~\eqref{eqn:model_problem} with $a
  = 10^{-p}$ and $b = 10^{p}$ for $p = 0, 1, \dotsc, 5$ from left to right. 
  The lower bounds at -50 are only active for $p = -2$ (deep blue).}
  \label{fig:optimal_controls}
\end{figure}

\begin{figure}[tb]
  \centering
  \includegraphics[width=\textwidth]{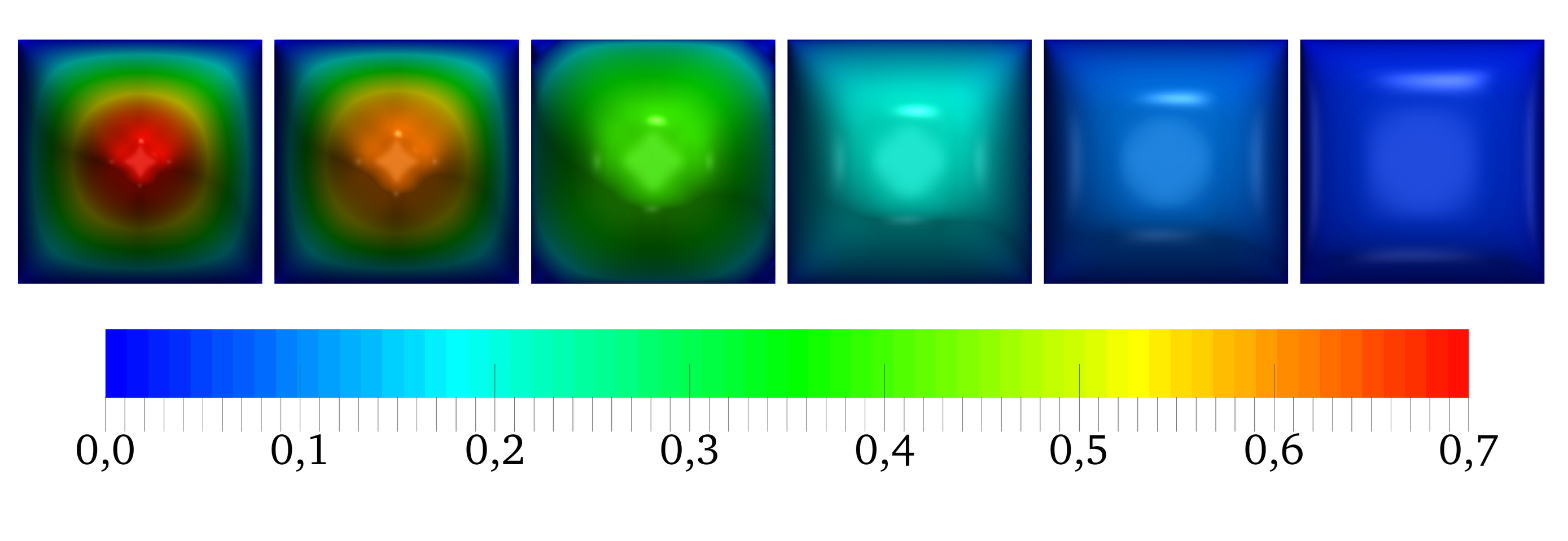}
  \caption{Optimal states for problem~\eqref{eqn:model_problem} with $a
  = 10^{-p}$ and $b = 10^{p}$ for $p = 0, 1, \dotsc, 5$ from left to right.}
  \label{fig:optimal_states}
\end{figure}

We apply Algorithm~\ref{alg:seqhom} to problem~\eqref{eqn:model_problem} on
$\Omega = (0,1)^2$ with the target state $u_{\mathrm{d}}(\xi) = 12
(1-\xi_1) \xi_1 (1-\xi_2) \xi_2$ from~\cite{Lubkoll2017} and control bounds
\begin{align*}
  q_{\mathrm{l}}(\xi) &= -50, &
  q_{\mathrm{u}}(\xi) &= \min \left (50, 800 \max \left( \left(
  \xi_1-\tfrac{1}{2} \right)^2, \left( \xi_2-\tfrac{1}{2} \right)^2 \right)
  \right)
\end{align*}
for the parameters $a = 10^{-p}$, $b = 10^{p}$ for $p = 0, \dotsc, 5$ with
continuous piecewise linear (P1) finite elements on regular triangular grids
with $N = 64, 128, 256, 512$ elements along each side of the unit square.

We perform Algorithm~\ref{alg:seqhom} with the initial guess $z_0 = 0$ and the
parameters $\Theta = 0.9$, $\lambda_{\mathrm{term}} = 10^{-8}$,
$\lambda_{\mathrm{inc}} = 2$, and $\mathrm{TOL} = 10^{-8}$. We fix
the choice of the penalty parameter to $\rho = 0.1$. For the stepsize PI
controller, we set $\theta_{\mathrm{ref}} = 0.5$, $K_{P} = 0.2$, $K_{I} =
0.005$, and $\lambda_{\mathrm{min}} = 10^{-12}$.
Figures~\ref{fig:optimal_controls} and~\ref{fig:optimal_states} depict
the resulting optimal controls and states.

\begin{table}[p]
  \centering
  \caption{Comparison of the sequential homotopy method of
    Algorithm~\ref{alg:seqhom} with a nonlinear VI solver with backtracking (bt)
  and error-oriented monotonicity test (nleqerr) for different instances of
  problem~\eqref{eqn:model_problem} and varying discretizations ($N$). The
  cardinality of the discrete optimal active set is given in the \#act column.
  The column \#disc shows the number of discarded steps, which are reasonably
  low, hinting at the efficiency of the PI control stepsize prediction.  The
  columns \#mat and \#res show the number of required matrix and residual
  evaluations.  The sequential homotopy method solves all instances and exhibits
  mesh-independent convergence (subject to some fluctuations for the worse
  conditioned problems). The symbol $\not\searrow$ denotes an error in the
  line-search, the symbol $\infty$ an error after exceeding 5.000 residual
  evaluations, and the symbol $\nearrow$ an error after not more than
  $10^{-8}$ relative reduction of a criticality measure over 100 system matrix
  evaluations.}
  \label{tab:control_constraints_decreasing_a}
\begin{tabular}{cccccccccccc}
\toprule
\multicolumn{3}{c}{Problem parameters} & Solution & \multicolumn{3}{c}{Sequential homotopy} & \multicolumn{2}{c}{VI (bt)} & \multicolumn{2}{c}{VI (nleqerr)}   \\
 $\log_{10} a$ &  $\log_{10} b$ &  $N$ &  \#act &  \#disc &  \#mat &  \#res &  \#mat &  \#res &  \#mat &  \#res \\
\midrule
%             0 &              0 &   32 &     157 &       0 &     20 &     40 &                9 &              26 &       \multicolumn{2}{c}{$\infty$} \\
             0 &              0 &   64 &     637 &       0 &     20 &     40 &                9 &              26 &     \multicolumn{2}{c}{$\nearrow$} \\
             0 &              0 &  128 &    2545 &       0 &     21 &     42 &                9 &              26 &     \multicolumn{2}{c}{$\nearrow$} \\
             0 &              0 &  256 &   10101 &       0 &     20 &     40 &                9 &              26 &               40 &             123 \\
             0 &              0 &  512 &   40193 &       0 &     20 &     40 &                9 &              26 &               12 &              33 \\
\midrule
%            -1 &              1 &   32 &     289 &       0 &     32 &     64 &               17 &              72 &     \multicolumn{2}{c}{$\nearrow$} \\
            -1 &              1 &   64 &    1121 &       0 &     32 &     64 &               16 &              67 &     \multicolumn{2}{c}{$\nearrow$} \\
            -1 &              1 &  128 &    4405 &       0 &     31 &     62 &               13 &              50 &     \multicolumn{2}{c}{$\nearrow$} \\
            -1 &              1 &  256 &   17525 &       0 &     32 &     64 &               12 &              42 &     \multicolumn{2}{c}{$\nearrow$} \\
            -1 &              1 &  512 &   69857 &       0 &     32 &     64 &               11 &              37 &               24 &              69 \\
\midrule
%            -2 &              2 &   32 &     749 &       3 &     55 &    113 &               44 &             332 &     \multicolumn{2}{c}{$\nearrow$} \\
            -2 &              2 &   64 &    2897 &       5 &     55 &    115 &     \multicolumn{2}{c}{$\nearrow$} &     \multicolumn{2}{c}{$\nearrow$} \\
            -2 &              2 &  128 &   11533 &      15 &     75 &    165 & \multicolumn{2}{c}{$\not\searrow$} &     \multicolumn{2}{c}{$\nearrow$} \\
            -2 &              2 &  256 &   45649 &       4 &     60 &    124 & \multicolumn{2}{c}{$\not\searrow$} &     \multicolumn{2}{c}{$\nearrow$} \\
            -2 &              2 &  512 &  182293 &       5 &     58 &    121 &     \multicolumn{2}{c}{$\nearrow$} &     \multicolumn{2}{c}{$\nearrow$} \\
\midrule
%            -3 &              3 &   32 &     905 &       1 &     46 &     93 &     \multicolumn{2}{c}{$\nearrow$} &              167 &             548 \\
            -3 &              3 &   64 &    3505 &       1 &     46 &     93 &     \multicolumn{2}{c}{$\nearrow$} &              110 &             378 \\
            -3 &              3 &  128 &   13997 &       1 &     47 &     95 &     \multicolumn{2}{c}{$\nearrow$} &     \multicolumn{2}{c}{$\nearrow$} \\
            -3 &              3 &  256 &   55709 &       4 &     55 &    114 &     \multicolumn{2}{c}{$\nearrow$} &     \multicolumn{2}{c}{$\nearrow$} \\
            -3 &              3 &  512 &  222385 &       3 &     54 &    111 &     \multicolumn{2}{c}{$\nearrow$} &     \multicolumn{2}{c}{$\nearrow$} \\
\midrule
%            -4 &              4 &   32 &     869 &       4 &     59 &    122 &     \multicolumn{2}{c}{$\nearrow$} &              105 &             369 \\
            -4 &              4 &   64 &    3405 &       4 &     59 &    122 &     \multicolumn{2}{c}{$\nearrow$} &       \multicolumn{2}{c}{$\infty$} \\
            -4 &              4 &  128 &   13477 &       4 &     56 &    116 &     \multicolumn{2}{c}{$\nearrow$} &               54 &             137 \\
            -4 &              4 &  256 &   53609 &       5 &     60 &    125 &     \multicolumn{2}{c}{$\nearrow$} &     \multicolumn{2}{c}{$\nearrow$} \\
            -4 &              4 &  512 &  214009 &       4 &     63 &    130 &     \multicolumn{2}{c}{$\nearrow$} &     \multicolumn{2}{c}{$\nearrow$} \\
\midrule
%            -5 &              5 &   32 &     749 &       7 &     72 &    151 &     \multicolumn{2}{c}{$\nearrow$} &               72 &             281 \\
            -5 &              5 &   64 &    2933 &      11 &     73 &    157 &     \multicolumn{2}{c}{$\nearrow$} &     \multicolumn{2}{c}{$\nearrow$} \\
            -5 &              5 &  128 &   11609 &      10 &     78 &    166 &     \multicolumn{2}{c}{$\nearrow$} &     \multicolumn{2}{c}{$\nearrow$} \\
            -5 &              5 &  256 &   46265 &      14 &     82 &    178 &     \multicolumn{2}{c}{$\nearrow$} &     \multicolumn{2}{c}{$\nearrow$} \\
            -5 &              5 &  512 &  184657 &      14 &     83 &    180 &     \multicolumn{2}{c}{$\nearrow$} &     \multicolumn{2}{c}{$\nearrow$} \\
\bottomrule
\end{tabular}
\end{table}

% Linear system condition number for a=1e-2, b=1e2, N=64: 5.01665e+12

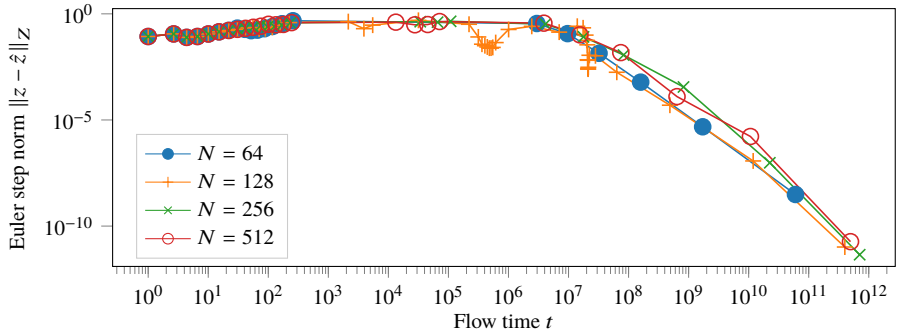
\begin{figure}[p]
  \centering
  % This file was created by matplotlib2tikz v0.6.18.
\begin{tikzpicture}

\definecolor{color0}{rgb}{0.12156862745098,0.466666666666667,0.705882352941177}
\definecolor{color1}{rgb}{1,0.498039215686275,0.0549019607843137}
\definecolor{color2}{rgb}{0.172549019607843,0.627450980392157,0.172549019607843}
\definecolor{color3}{rgb}{0.83921568627451,0.152941176470588,0.156862745098039}

\begin{axis}[
legend cell align={left},
legend entries={{$N = 64$},{$N = 128$},{$N = 256$},{$N = 512$}},
legend style={at={(0.03,0.03)}, anchor=south west, draw=white!80.0!black},
tick align=outside,
tick pos=left,
width=\textwidth,
height=5cm,
x grid style={white!69.01960784313725!black},
xlabel={Flow time $t$},
xmin=0.255600078030223, xmax=2761972667950.49,
xmode=log,
y grid style={white!69.01960784313725!black},
ylabel={Euler step norm $\norm{z-\hat{z}}{Z}$},
ymin=1.25518403873493e-12, ymax=1.76370359953831,
ymode=log
]
\addlegendimage{mark=*, color0}
\addlegendimage{mark=+, color1}
\addlegendimage{mark=x, color2}
\addlegendimage{mark=o, color3}
\addplot [semithick, color0, mark=*, mark size=3, mark options={solid}]
table [row sep=\\]{%
1	0.0855287303458733 \\
2.6579201199794	0.113188388154832 \\
4.37339521906494	0.0778684407859206 \\
6.61212942364445	0.0856049615379028 \\
10.0419243283172	0.112476329134311 \\
15.1375532590791	0.142117077617343 \\
21.9541450656086	0.165995905828095 \\
30.5850397878661	0.212474314637018 \\
40.7540478879046	0.186289973269715 \\
51.6807044662435	0.161194422270869 \\
66.1021410803125	0.170322627837922 \\
86.6775616467883	0.199359573677192 \\
118.219187222762	0.254550220904914 \\
169.524265516541	0.34212403271232 \\
257.493964229717	0.47207829551971 \\
2968773.33024447	0.350236802745107 \\
9657429.37264612	0.116899549003409 \\
32186026.6596586	0.013827373394592 \\
159073901.830059	0.000605155200765424 \\
1715602664.31916	4.7681936345294e-06 \\
60301677067.9314	3.07806744669915e-09 \\
};
\addplot [semithick, color1, mark=+, mark size=3, mark options={solid}]
table [row sep=\\]{%
1	0.085544290582323 \\
2.65759628281059	0.113163702592613 \\
4.37260848690311	0.0778413444633828 \\
6.61077726905513	0.0854438017106832 \\
10.0397880514007	0.112298624065784 \\
15.1373512326169	0.141777895280178 \\
21.9710780775418	0.164750717661326 \\
30.6686943956966	0.186854957544646 \\
41.7037493546585	0.217912092691055 \\
55.545394255302	0.226218079472536 \\
70.3311693498935	0.191809170048216 \\
89.4620027432287	0.19889650378889 \\
116.579348691752	0.229875467483466 \\
157.517749508055	0.29020629044071 \\
223.960896277027	0.387555579174359 \\
2151.10016056858	0.430222844457815 \\
3914.59105981607	0.200218066606385 \\
5481.50997384513	0.293116836171215 \\
31644.6561953428	0.49460282489056 \\
220352.19287671	0.323075136248522 \\
314705.961217394	0.081582143415983 \\
361882.845387736	0.0371247250870291 \\
404882.779734517	0.0315054737881496 \\
443880.499833079	0.0274932724495433 \\
479376.580083405	0.0242407089675129 \\
512478.442256523	0.0243546257013569 \\
547943.277934222	0.0261971391870167 \\
612604.470269924	0.0442455855652254 \\
1005109.23488739	0.183442042546097 \\
2443397.3895129	0.258174046980096 \\
6991378.82795753	0.137530788532209 \\
13995596.1341997	0.259798676381184 \\
17497704.7873208	0.219236380053249 \\
19248759.1138813	0.0994773620254868 \\
20124286.2771616	0.0563969591149796 \\
20562049.8588017	0.0346918495545869 \\
20780931.6496218	0.00685957148539949 \\
20977764.7854942	0.00650740495614701 \\
21076181.3534304	0.00266656381015436 \\
21164974.0973986	0.00257819514312867 \\
21256126.5220805	0.00248452077482108 \\
21383288.4176071	0.00304275171327593 \\
22285703.8948632	0.0111721607384212 \\
28071989.630885	0.0105246300399356 \\
63661033.507071	0.0017671248156072 \\
487740952.788596	4.95510940011356e-05 \\
11887309353.9617	1.17109625552365e-07 \\
399643777045.576	1.0315657028799e-11 \\
};
\addplot [semithick, color2, mark=x, mark size=3, mark options={solid}]
table [row sep=\\]{%
1	0.0855481810858469 \\
2.65751652819619	0.113163824722907 \\
4.37238132622381	0.0778383414489905 \\
6.61038034841672	0.085549614708702 \\
10.0404119581934	0.11242381520606 \\
15.1405514615391	0.141906708132468 \\
21.9804575266669	0.164795196554849 \\
30.7022753064393	0.185904963818516 \\
41.8242427145423	0.208397014594271 \\
56.2295037614978	0.234859338097508 \\
75.0220255428482	0.284767961758104 \\
97.19628571019	0.243392915483119 \\
126.124935723951	0.250458770979922 \\
167.626921715438	0.293991214816009 \\
231.944701896541	0.37804846237339 \\
33929.0102519764	0.426298743558665 \\
65974.8496537533	0.405346046416121 \\
108687.590506717	0.444929268061015 \\
4111779.90836871	0.328077181779991 \\
17766151.4813139	0.0820400069548058 \\
82593408.823512	0.0114509646438305 \\
822280645.642071	0.000357740266859504 \\
22504638454.1544	9.59132500048514e-08 \\
705960429445.492	4.47585920620177e-12 \\
};
\addplot [semithick, color3, mark=o, mark size=3, mark options={solid}]
table [row sep=\\]{%
1	0.0855491536754131 \\
2.65749635143522	0.113162596159891 \\
4.37233064423593	0.0778381155637726 \\
6.61029487467852	0.0855435611532459 \\
10.0403494743633	0.112397766134603 \\
15.1407289250841	0.14182753990752 \\
21.981913870281	0.164659773563484 \\
30.7075711896829	0.185628230194651 \\
41.8465878048414	0.207596788548935 \\
56.3241610887621	0.233075121479305 \\
75.4391381150431	0.261333188405479 \\
100.951257918536	0.322317758793752 \\
134.213813398366	0.304847103133576 \\
178.172919352127	0.316136590504302 \\
244.027659935438	0.385014457349213 \\
13304.6381580357	0.405775884279906 \\
27621.9063145363	0.304563914301885 \\
45184.0466898621	0.313270189456558 \\
71250.7912800302	0.433923767530527 \\
3891031.96131557	0.363632066919146 \\
15487418.8049286	0.10227895877625 \\
74257051.0180812	0.0149065860457665 \\
640908628.89311	0.000125998993923165 \\
10700797750.3029	1.66497224096203e-06 \\
496803400100.295	1.8465184359828e-11 \\
};
\end{axis}

\end{tikzpicture}
  \caption{The projected backward Euler step norms $\norm{z - \hat{z}}{Z}$
  for~\eqref{eqn:model_problem} with $a = 10^{-2}$, $b = 10^{2}$ on
  different meshes plotted with respect to the flow time $t$, which is the sum
  of all accepted step sizes $\Delta t = 1/\lambda$.}
  \label{fig:convergence}
\end{figure}

%\afterpage{\clearpage}

We compare the sequential homotopy method of Algorithm~\ref{alg:seqhom} with a
nonlinear VI solver described in~\cite{Munson2001,Benson2006} and implemented in
the production quality software package
PETSc~\cite{petsc-web-page,petsc-user-ref}. For better comparison, we use the
direct solver MUMPS~\cite{Amestoy2001,Amestoy2006} for the solution of the
linear systems in both approaches.
The use of inexact linear algebra solvers is no conceptual problem, as long as
they yield a locally convergent nonlinear iteration. The efficiency of
iterative linear algebra methods, however, depends crucially on the use of
suitable structure-exploiting preconditioners. This topic exceeds the
scope of this paper and is the subject of future research.

For the VI solver, we consider two implemented globalization strategies, a
backtracking line-search (bt) and an error-oriented monotonicity test (nleqerr).
As it turns out, the VI solver did not solve any of the problem instances when
started at the initial guess $z_0 = 0$, failing
either by raising an error or reaching the limit of 5.000 residual
evaluations, even for a reduced termination tolerance of $10^{-5}$ on the
$l^{\infty}$-norm of the residuals. Some problem instances could be solved
successfully after dropping the lower control bound, which is only active for
$a = 10^{-2}$, $b = 10^{2}$. In some of these instances the residual norm
stalled between $10^{-5}$ and $10^{-8}$.

We compare in Table~\ref{tab:control_constraints_decreasing_a} the sequential
homotopy method of Algorithm~\ref{alg:seqhom} (with a
sharper termination tolerance of $10^{-8}$ on the $Z$-norm of the homotopy
increment and upper and lower bounds) to the VI approach with reduced
termination tolerance as above and only upper bounds. We can observe that the
sequential homotopy method solves all problem instances with mesh-independent
convergence (subject to some fluctuation for the worse conditioned problems).
The VI approach with backtracking is faster for the less demanding but fails for
the more demanding instances. The VI approach with error-oriented monotonicity
test solves at least two of the more demanding instances successfully, although
only one with an efficiency comparable to the sequential homotopy method.

In Figure~\ref{fig:convergence}, we see that even though slightly different
numbers of iterations (depicted with markers) are performed on different meshes
for the case $a = 10^{-2}$, $b=10^{2}$, roughly the same flow time of $10^{11}$
has to be traversed to reach the required tolerance of $\mathrm{TOL}=10^{-8}$.
We also see that the stepsizes $\Delta t$ eventually become very large and lead
to superlinear convergence. This is the typical numerical behavior of the
sequential homotopy method on all considered instances. For $N=128$ some extra
steps are carried out around $t=10^{5}$ and $t=10^{7}$.

\section{Summary}

We provided sufficient conditions for the existence of global solutions to the
projected gradient/antigradient flow~\eqref{eqn:gradflow} and showed that
critical points
with emanating descent curves cannot be asymptotically stable and are thus not
attracting for the flow. We applied projected backward Euler timestepping to
derive the necessary optimality conditions of a primal-dual proximally
regularized counterpart~\eqref{eqn:homotopy_mp} of~\eqref{eqn:mp}. The
regularized problem can be solved by a homotopy method, giving rise to a
sequence of homotopy problems. The sequential homotopy method can be used to
globalize any locally convergent optimization method that can be employed
efficiently in a homotopy framework. The sequential homotopy method with a local
semismooth Newton solver outperforms state-of-the-art VI solvers for a
challenging class of PDE-constrained optimization problem with control
constraints.

%\begin{figure}
%  \includegraphics{example.eps}
%\end{figure}
%\begin{table}
%\begin{tabular}{lll}
%\hline\noalign{\smallskip}
%\end{tabular}
%\end{table}

%\begin{acknowledgements}
%\end{acknowledgements}

%\bibliography{seqhom,fenics}   % name your BibTeX data base
%#\input{seqhom.bbl}

\begin{center}
  \Large
  Correction to:\\ Mathematical Programming (2021) 187:459–486 https://doi.org/10.1007/s10107-020-01488-z
\end{center}

\paragraph{Correction 1}
For the example~(26) considered in Sec.~5 of ``A sequential homotopy method for mathematical programming problems'' written by Potschka, A.~and Bock, H.G., the argument presented in the last paragraph of Sec.~5.1 is insufficient to prove semismoothness of the equations in~(21). The required smoothing property for the argument of
$
  P_{C_Q}\left(\hat{q} - \Delta t \left[ \gamma q - R_U (y + \rho c(x))\right]\right)
$
cannot be established with the usual trick (see, e.g., \cite[Theorem~2.14]{hinze2009optimization}) of letting $\Delta t = 1 / \gamma$, which only works if the first term in the argument is also $q$ (and not $\hat{q}$). In fact, semismoothness does not hold due to~\cite[Lemma~2.7]{hinze2009optimization}.

Hence, semismoothness of~(21) holds only for each fixed discretization and there is no theoretical justification for the good numerical results for refined discretizations reported in Sec.~5.3 using a semismooth Newton method on the projected backward Euler equations~(21) for the gradient/antigradient flow of the augmented Lagrangian directly.

This gap can be closed with an algorithmic modification to solve the homotopy subproblem~(22) differently (i.e., not by (21) directly): By Lemma~5, critical points of~(22) can also be equivalently characterized as equilibria of another gradient/antigradient flow for~(22), which are equivalent with the fixpoints of the projected (forward or backward) Euler equations for any fixed stepsize $\tau > 0$. Using the Lagrangian $L^{\lambda,\rho}(w, x, y)$ of (22), which is defined at the beginning of the proof of Theorem~4, the projected Euler fixpoint equations for~(22) read
\begin{align}
  w &= w - \tau \nabla_w L^{\rho,\lambda}(w, x, y), \label{eqn:euler_fixpoint_w}\tag{C1a} \\ % &&\Leftrightarrow & 0 &= w + y - \hat{y} \\
  x &= P_C\left(x - \tau \nabla_x L^{\rho,\lambda}(w, x, y)\right), \label{eqn:proj_x} \tag{C1b} \\ % &&\Leftrightarrow & 0 &= x - P_C\left((1-\tau \lambda) x - \tau \lambda \hat{x} - \theta \nabla_x L^{\rho}(x,y)\right) \\
  y &= y + \tau \nabla_y L^{\rho,\lambda}(w, x, y). \label{eqn:euler_fixpoint_y}\tag{C1c} % \tau (c(x) + \lambda w). % &&\Leftrightarrow & 0 &= c(x) + \lambda w
\end{align}
For the problem class~(26), the projection acts only on the control component $q$ of $x = (u, q)$ and
\[
  \nabla_q L^{\rho,\lambda}(w, x, y) = \gamma q - R_U (y + \rho c(x)) + \lambda (q - \hat{q}).
\]
Hence, the corresponding line in~\eqref{eqn:proj_x} expands to 
\begin{align*}
  q &= P_{C_Q} \left(q - \tau [\gamma q - R_U (y + \rho c(x)) + \lambda (q - \hat{q})] \right) \\
  &= P_{C_Q} \left( [1 - \tau (\gamma + \lambda)] q + \tau \lambda \hat{q} + \tau R_U (y + \rho c(x))\right).
\end{align*}
Choosing $\tau = \frac{1}{\gamma + \lambda}$ cancels the $q$-term and if $C_Q \subset L^r(\Omega)$ for some $r > 2$, the required smoothing property holds because $\hat{q} \in C_Q$ and the Riesz operator $R_U$ maps to $H^1_0(\Omega)$.

A semismooth Newton method can then be applied to~\eqref{eqn:euler_fixpoint_w}--\eqref{eqn:euler_fixpoint_y} instead of~(21), where~\eqref{eqn:euler_fixpoint_w} can be used to immediately eliminate $w = \hat{y} - y$ just as in~(24). Using
\begin{align*}
  \nabla_x L^{\rho,\lambda}(w, x, y) &= \nabla_x L^{\rho}(x, y) + \lambda (x - \hat{x}), \\
  \nabla_y L^{\rho,\lambda}(w, x, y) &= c(x) + \lambda w = c(x) - \lambda (y - \hat{y}),
\end{align*}
and scaling by $1/\tau$ delivers the system of equations
\begin{align}
  0 &= \tfrac{1}{\tau} x - \tfrac{1}{\tau} P_C((1 - \tau \lambda) x + \tau \lambda \hat{x} - \tau \nabla_x L^{\rho}(x, y)), \tag{C2a} \label{eqn:ssn_x} \\
  0 &= c(x) - \lambda (y - \hat{y}). \tag{C2b} \label{eqn:ssn_y}
\end{align}
The resulting linear subproblems of a semismooth Newton method for~\eqref{eqn:ssn_x}--\eqref{eqn:ssn_y} (and, equivalently, \eqref{eqn:euler_fixpoint_w}--\eqref{eqn:euler_fixpoint_y}), differ from the ones given in Sec.~5.1 only in the determination of the active set, while the remaining entries of the matrices and right-hand sides coincide. The original formulation~(21) can be recovered by choosing $\tau = \Delta t = \frac{1}{\lambda}$ instead of $\tau = \frac{1}{\gamma + \lambda}$, which justifies good numerical behavior for $\lambda \gg \gamma$.

We provide an update of Table~\ref{tab:results} with the results for the modified active set determination. In terms of number of matrix evaluations, the results are worse for $b=10^0$, similar for $b=10^1$, slightly better for $b=10^2$, slightly worse for $b=10^3$ (with an outlier on the $N=512$ mesh, which vanishes for a small perturbation of the initial $\lambda=1.1$ instead of $\lambda=1$), slightly better for $b=10^4$, and considerably better for $b=10^5$. Except for the outlier, the number of required matrix evaluations appears to be mesh independent.

We believe the reason for the outlier to be the following: It may happen that the monotonicity test accepts a step even though the inertia of the saddle-point matrix of the linearized subproblem for (21) changes between the points $z$, $z^+$ and $z^{++}$, which implies the existence of a singularity on one of the lines connecting $z$ with $z^+$ or $z^+$ with $z^{++}$. As a result, $\lambda$ needs to be increased considerably, resulting in a high number of discarded steps. It would generally be possible to check the inertia using appropriate sparse matrix decomposition methods, but it is technically challenging in our current implementation. It is reassuring that the sequential homotopy method eventually recovers and converges to a solution (in this case the same as in the original version).

\begin{table}[htbp]
  \begin{tabular}{cccccccccccc}
    \toprule
    \multicolumn{3}{c}{Problem parameters} & Solution & \multicolumn{3}{c}{Sequential homotopy} & \multicolumn{2}{c}{VI (bt)} & \multicolumn{2}{c}{VI (nleqerr)}   \\
    $\log_{10} a$ &  $\log_{10} b$ &  $N$ &  \#act &  \#disc &  \#mat &  \#res &  \#mat &  \#res &  \#mat &  \#res \\
    \midrule
                0 &              0 &   64 &     637 &       0 &\st{20} 24& \st{40}  47&                9 &              26 &     \multicolumn{2}{c}{$\nearrow$} \\
                0 &              0 &  128 &    2545 &       0 &\st{21} 25& \st{42}  49&                9 &              26 &     \multicolumn{2}{c}{$\nearrow$} \\
                0 &              0 &  256 &   10101 &       0 &\st{20} 26& \st{40}  51&                9 &              26 &               40 &             123 \\
                0 &              0 &  512 &   40193 &       0 &\st{20} 26& \st{40}  51&                9 &              26 &               12 &              33 \\ \midrule
                -1 &              1 &   64 &    1121 &       0 &\st{32} 31& \st{64}  61&               16 &              67 &     \multicolumn{2}{c}{$\nearrow$} \\
                -1 &              1 &  128 &    4405 &       0 &\st{31} 32& \st{62}  63&               13 &              50 &     \multicolumn{2}{c}{$\nearrow$} \\
                -1 &              1 &  256 &   17525 &       0 &\st{32} 32& \st{64}  63&               12 &              42 &     \multicolumn{2}{c}{$\nearrow$} \\
                -1 &              1 &  512 &   69857 &       0 &\st{32} 32& \st{64}  63&               11 &              37 &               24 &              69 \\ \midrule
                -2 &              2 &   64 &    2897 &\st{ 5} 2&\st{55} 52&\st{115} 105&     \multicolumn{2}{c}{$\nearrow$} &     \multicolumn{2}{c}{$\nearrow$} \\
                -2 &              2 &  128 &   11533 &\st{15} 1&\st{75} 51&\st{165} 102& \multicolumn{2}{c}{$\not\searrow$} &     \multicolumn{2}{c}{$\nearrow$} \\
                -2 &              2 &  256 &   45649 &\st{ 4} 7&\st{60} 59&\st{124} 124& \multicolumn{2}{c}{$\not\searrow$} &     \multicolumn{2}{c}{$\nearrow$} \\
                -2 &              2 &  512 &  182293 &\st{ 5} 6&\st{58} 57&\st{121} 119&     \multicolumn{2}{c}{$\nearrow$} &     \multicolumn{2}{c}{$\nearrow$} \\ \midrule
                -3 &              3 &   64 &    3505 & \st{1} 2&\st{46} 48&\st{ 93}  97&     \multicolumn{2}{c}{$\nearrow$} &              110 &             378 \\
                -3 &              3 &  128 &   13997 & \st{1} 1&\st{47} 51&\st{ 95} 102&     \multicolumn{2}{c}{$\nearrow$} &     \multicolumn{2}{c}{$\nearrow$} \\
                -3 &              3 &  256 &   55709 & \st{4} 1&\st{55} 53&\st{114} 106&     \multicolumn{2}{c}{$\nearrow$} &     \multicolumn{2}{c}{$\nearrow$} \\
                -3 &              3 &  512 &  222385 & \st{3} 94&\st{54} 241&\st{111} 576&     \multicolumn{2}{c}{$\nearrow$} &     \multicolumn{2}{c}{$\nearrow$} \\
                -3 &              3 &  512 &  222385 & 1$^\ast$& 53$^\ast$& 106$^\ast$ \\ \midrule
                -4 &              4 &   64 &    3405 & \st{4} 2&\st{59} 51&\st{122} 103&     \multicolumn{2}{c}{$\nearrow$} &       \multicolumn{2}{c}{$\infty$} \\
                -4 &              4 &  128 &   13477 & \st{4} 5&\st{56} 59&\st{116} 122&     \multicolumn{2}{c}{$\nearrow$} &               54 &             137 \\
                -4 &              4 &  256 &   53609 & \st{5} 3&\st{60} 56&\st{125} 114&     \multicolumn{2}{c}{$\nearrow$} &     \multicolumn{2}{c}{$\nearrow$} \\
                -4 &              4 &  512 &  214009 & \st{4} 2&\st{63} 56&\st{130} 113&     \multicolumn{2}{c}{$\nearrow$} &     \multicolumn{2}{c}{$\nearrow$} \\ \midrule
                -5 &              5 &   64 &    2933 &\st{11} 5&\st{73} 59&\st{157} 122&     \multicolumn{2}{c}{$\nearrow$} &     \multicolumn{2}{c}{$\nearrow$} \\
                -5 &              5 &  128 &   11609 &\st{10} 6&\st{78} 66&\st{166} 137&     \multicolumn{2}{c}{$\nearrow$} &     \multicolumn{2}{c}{$\nearrow$} \\
                -5 &              5 &  256 &   46265 &\st{14} 5&\st{82} 64&\st{178} 132&     \multicolumn{2}{c}{$\nearrow$} &     \multicolumn{2}{c}{$\nearrow$} \\
                -5 &              5 &  512 &  184657 &\st{14} 5&\st{83} 66&\st{180} 136&     \multicolumn{2}{c}{$\nearrow$} &     \multicolumn{2}{c}{$\nearrow$} \\
    \bottomrule
  \end{tabular}
  \caption{Updated Table~1. Run with an initial $\lambda=1.1$ instead of $\lambda=1$ marked with $^\ast$.}
  \label{tab:results}
\end{table}

%  0 &  0 &   64 &    637 &  0 &  24 &  47
%  0 &  0 &  128 &   2545 &  0 &  25 &  49
%  0 &  0 &  256 &  10101 &  0 &  26 &  51
%  0 &  0 &  512 &  40193 &  0 &  26 &  51
% -1 &  1 &   64 &   1121 &  0 &  31 &  61
% -1 &  1 &  128 &   4405 &  0 &  32 &  63
% -1 &  1 &  256 &  17525 &  0 &  32 &  63
% -1 &  1 &  512 &  69857 &  0 &  32 &  63
% -2 &  2 &   64 &   2897 &  2 &  52 & 105
% -2 &  2 &  128 &  11533 &  1 &  51 & 102
% -2 &  2 &  256 &  45649 &  7 &  59 & 124
% -2 &  2 &  512 & 182293 &  6 &  57 & 119
% -3 &  3 &   64 &   3505 &  2 &  48 &  97
% -3 &  3 &  128 &  13997 &  1 &  51 & 102
% -3 &  3 &  256 &  55709 &  1 &  53 & 106
% -3 &  3 &  512 & 222385 &  1 &  53 & 106 % lambda_0 = 1.1
% -3 &  3 &  512 & 222385 & 94 & 241 & 576
% -4 &  4 &   64 &   3405 &  2 &  51 & 103
% -4 &  4 &  128 &  13477 &  5 &  59 & 122
% -4 &  4 &  256 &  53609 &  3 &  56 & 114
% -4 &  4 &  512 & 214009 &  2 &  56 & 113
% -5 &  5 &   64 &   2933 &  5 &  59 & 122
% -5 &  5 &  128 &  11609 &  6 &  66 & 137
% -5 &  5 &  256 &  46265 &  5 &  64 & 132
% -5 &  5 &  512 & 184657 &  5 &  66 & 136

\paragraph{Correction 2} Substitute ``directional'' for ``Gâteaux'' in Lemma~6, because the limit is only taken over positive $h$.

% \bibliographystyle{plain}
% \bibliography{refs}

\end{document}